\input amstex
\documentstyle{amsppt}
\NoBlackBoxes

\TagsOnRight

\def\cal{\Cal}
\def\AA{{\cal A}}

\def\HH{{\cal H}}
\def\MM{{\cal M}}

\def\JJ{{\cal J}}
\def\UU{{\cal U}}

\def\UU{{\cal U}}

\def\Z{{\Bbb Z}}
\def\C{{\Bbb C}}
\def\R{{\Bbb R}}

\def\e{{\epsilon}}

\def\n{\noindent}
\def\part{{\partial}}
\def\dudtau{{\part u\over \part \tau}}
\def\dudt{{\part u\over \part t}}

\rightheadtext{Perturbed Cauchy-Riemann equation} \leftheadtext{
Yong-Geun Oh }

\topmatter
\title
An existence theorem, with energy bounds, of Floer's perturbed
Cauchy-Riemann equation with jumping discontinuity
\endtitle
\author
Yong-Geun Oh\footnote{Partially supported by the NSF Grant \#
DMS-9971446 \& DMS-0203593 and by \# DMS-9729992 in the Institute for Advanced
Study, Vilas Associate Award in the University of Wisconsin
and by a grant of the Korean Young Scientist Prize
\hskip6.5cm\hfill}
\endauthor
\address
Department of Mathematics, University of Wisconsin, Madison, WI
53706, ~USA \& Korea Institute for Advanced Study, Seoul, Korea;
oh\@math.wisc.edu
\endaddress
\abstract This is a sequel to the paper [Oh5]. The main purpose of
the present paper is to give the proof of an existence theorem
{\it with energy bounds} of certain pseudo-holomorphic sections of
the deformed mapping cylinder that is needed for the proof of
nondegeneracy of the homological invariant pseudo-norm which the
author has constructed on general symplectic manifolds [Oh4,5].
The existence theorem is also the crux of the author's recent
proof of an optimal energy-capacity inequality given in [Oh5]. In
this paper, we give a more general existence result than needed in
that we study Floer's perturbed Cauchy-Riemann equations with {\it
discontinuous} Hamiltonian perturbation terms and prove an
existence theorem of certain piecewise smooth finite energy
solutions of the equation. The proof relies on a careful study of
the product structure on the Floer homology and a singular
degeneration ({\it ``adiabatic degeneration''}) of Floer's
perturbed Cauchy-Riemann equation. In the course of the proof, we
also derive certain general energy identity of pseudo-holomorphic
sections of the Hamiltonian fibration. 
\endabstract


\endtopmatter

\document

\bigskip

\centerline{\bf Contents} \medskip

\n \S1. Introduction and the main results
\smallskip

\n \S2. Perturbed Cauchy-Riemann equations re-visited
\smallskip

\n \S3. Energy identity on the deformed mapping cylinder
\smallskip

\n \S4. Pants product and the Hamiltonian fibration
\smallskip

\n \S5. Construction of the $W^{1,2}$-sections I: analysis of the
thick part
\smallskip

\n \S6. Construction of the $W^{1,2}$-sections II: analysis of the
thin part
\smallskip

\n \S7. The case with $H_3 = H_1\# H_2$
\smallskip

\vskip0.5in \head {\bf \S 1. Introduction}
\endhead

In [Oh3,4], on general (non-exact) compact symplectic manifolds,
the present author defined spectral invariants
$$
\rho: \widetilde{Ham}(M,\omega) \times QH^*(M)\setminus \{0\} \to
\R \tag 1.1
$$
by constructing a function
$$
\rho_a = \rho(\cdot; a): C^0_m(S^1 \times M)  \to \R
$$
for each $a \in QH^*(M) \setminus \{0\}$ such that
$$
\rho(H;a) \in \text{Spec}(H)
$$
and
$$
\rho(H;a) = \rho(F;a)
$$
if $[H] = [F]$, when $H, \, F$ are $C^2$-functions. The map (1.1)
is then defined by
$$
\rho(h;a): = \rho(H;a)
$$
for $h \in \widetilde{\HH am}(M,\omega)$ and $[H] = h$. We also
prove that $\rho_a$ is $C^0$-continuous and the invariant (1.1)
satisfies certain axioms which we refer to [Oh4].

Among these invariants, the invariant $\rho(h;1)$ or $\rho(H;1)$
is of particular interest. Using them, we [Oh5] constructed an
invariant pseudo-norm by considering the sum
$$
\widetilde \gamma(h) = \rho(h;1) + \rho(h^{-1};1) \tag 1.2
$$
and taking the infimum
$$
\gamma(\phi) = \inf_{\pi(h) = \phi} \widetilde \gamma(h). \tag 1.3
$$
The  inequality
$$
\aligned
\rho(h;1) & \leq E^-(h) : =  \inf_{[K] = h} E^-(K) \\
\rho(h^{-1};1) & \leq E^+(h) : = \inf_{[K] = h} E^+(K)
\endaligned
\tag 1.4
$$
is also proved in [Oh5] in general. Here we denote
$$
E^+(K) = \int_0^1 \max K_t\, dt, \quad
E^-(K) = \int_0^1 -\min K_t\,dt.
$$
In particular we have
$$
\rho(H;1) \leq E^-(h),\quad \rho(\widetilde F;1) \leq E^+(f)
\tag 1.5
$$
for any $H, \, F$ with $[H] = h$, $[F]=f$.

One crucial ingredient in the nondegeneracy proof of $\gamma: \HH
am(M,\omega) \to \R_+$ given in [Oh5] is the usage of a geometric
invariant $A(\phi,J_0)$ of the pair $(\phi,J_0)$ and certain
existence theorem of the perturbed Cauchy-Riemann equation for
nondegenerate Hamiltonian $H$. One of the purposes of the present
paper is to provide a proof of this existence theorem.

We first recall the definition of the invariant $A(\phi,J_0)$ from
[Oh5]. Denote by $J_0$ a compatible almost complex structure on
$(M,\omega)$ and by $\JJ_\omega$ the set of compatible almost
complex structures on $M$. For given $\phi$ and $J_0$, we consider
the set of paths $J': [0,1] \to \JJ_\omega$ with
$$
J'(0) = J_0, \quad J'(1) = \phi^*J_0  \tag 1.6
$$
and denote the set of such paths by
$$
j_{(\phi,J_0)}.
$$

For each given $J' \in j_{(\phi,J_0)}$, we define the constant
$$
\aligned A_S(\phi,J_0;J') = \inf \{\omega([u]) & \mid  u: S^2 \to
M \text{
non-constant and} \\
& \text{ satisfying $\overline \part_{J'_t}u = 0$ for some $t \in
[0,1]$}\}
\endaligned
\tag 1.7
$$
and then
$$
A_S(\phi,J_0) = \sup_{J'\in j_{(\phi,J_0)}} A_S(\phi,J_0;J'). \tag
1.8
$$
As usual, we set $A_S(\phi,J_0) =\infty$ if there is $J \in
j_{(\phi,J_0)}$ for which there is no $J_t$-holomorphic sphere for
any $t \in [0,1]$ as in the weakly exact case. The positivity
$A_S(\phi,J_0;J') > 0$ {\it when it is not infinite} and so
$A_S(\phi,J_0)> 0$ is an immediate consequence of the one
parameter version of the uniform $\e$-regularity theorem (see
[SU], [Oh1]).

Next for each given $J'\in j_{(\phi,J_0)}$, we consider the
equation of $v: \R \times [0,1] \to M$
$$
\cases {\part v \over \part \tau} + J'_t {\part v \over
\part t} = 0 \\
\phi(v(\tau,1)) = v(\tau,0), \quad \int |{\part v \over \part \tau
}|_{J'_t}^2 < \infty.
\endcases
\tag 1.9
$$
This equation itself is analytically well-posed and (1.6) enables
us to interpret solutions of (1.9) as pseudo-holomorphic sections
of the mapping cylinder of $\phi$ with respect to a suitably
chosen almost complex structure $\widetilde J$ on the mapping cylinder.

Note that any such solution of (1.9) also has the limit
$\lim_{\tau \to \pm}v(\tau)$ and satisfies
$$
\lim_{\tau \to \pm\infty}v(\tau) \in \text{Fix}\phi.
$$
Now it is a crucial matter to produce non-constant solutions of
(1.9), when $\phi$ is not the identity and in particular when
$\phi$ is a nondegenerate Hamiltonian diffeomorphism.

Suppose that $\phi \neq id$ is nondegenerate and choose a
symplectic ball $B(u)$ such that
$$
\phi(B(u)) \cap B(u) = \emptyset \tag 1.10
$$
where $B(u)$ is the image of a symplectic embedding into $M$ of
the standard Euclidean ball of radius $r$ with $u = \pi r^2$.
We then study (1.9) together with the condition
$$
v(0,0) \in B(u). \tag 1.11
$$
Because of (1.10), it follows
$$
v(\pm\infty) \in \text{Fix }\phi \subset M \setminus B(u). \tag
1.12
$$
Therefore such a solution cannot be constant because of (1.11) and
(1.12).

We now define the  constant
$$
A_D(\phi,J_0;J'): = \inf_{v} \Big\{ \int v^*\omega,  \mid
\text{$v$ non-constant solution of (1.9)}\Big\} \tag 1.13
$$
for each $J' \in j_{(\phi,J_0)}$. Again we have $A_D(\phi,J_0;J')
> 0$. We set
$$
A(\phi,J_0;J') = \min\{A_S(\phi, J_0;J'), A_D(\phi,J_0;J')\} \tag
1.14
$$
The following existence theorem is an immediate consequence of
Theorem II below.

\proclaim{Theorem I} For any nondegenerate Hamiltonian
diffeomorphism $\phi \neq id$, we have
$$
0 < A(\phi,J_0;J') < \infty \tag 1.15
$$
for any $J_0\in \JJ_\omega$ and $J' \in j_{(\phi,J_0)}$. More
precisely, the following alternative holds: for any given point in
$B(u)$ with
$$
\phi(B(u)) \cap B(u) = \emptyset,
$$

$(1)$ either there exists a non-constant $J'_t$-holomorphic sphere
for some $t \in [0,1]$ that pass through the point or

$(2)$ (1.9) has a non-constant solution that pass through the
point.
\endproclaim

In fact, Theorem II was also the crux in the proof of the
following optimal energy-capacity inequality which in particular
gives rise to nondegeneracy of the pseudo-norm $\gamma$. We refer
to [Oh5] for its proof based on Theorem II.

\proclaim{[Optimal Energy-Capacity inequality] [Oh5]} We denote by
$e_\gamma(A)$ the {\it $\gamma$-displacement energy} of closed
subset $A$ of $M$, i.e.,
$$
e_\gamma(A):= \inf\{ \gamma(\phi) \mid \phi(A) \cap A = \emptyset
\}
$$
and by $c(A)$ the Gromov capacity
$$
c(A):= \sup\{u\mid \exists \, \text{\rm a symplectic embedding}\, B(u)
\hookrightarrow \operatorname{Int}A \}.
$$
Then we have following inequalities:
$$
e_\gamma(A) \geq c(A)
$$
\endproclaim
As we pointed out in [Oh5], this gives rise to the optimal form of
the inequality between the Hofer displacement energy
and the Gromov capacity
$$
e(A) \geq c(A)
$$
since we have $e(A) \geq e_\gamma(A)$, where
$$
e(A) = \inf_{\phi}\{\|\phi\|\mid \phi(A)\cap A =\emptyset\}.
$$

Finally we define
$$
A(\phi,J_0) : = \sup_{J' \in j_{(\phi,J_0)}} \min\{A_S(\phi,
J_0;J'), A_D(\phi,J_0;J')\} \tag 1.16
$$
and
$$
A(\phi) = \sup_{J_0} A(\phi,J_0).
$$
Because of the assumption that $\phi$ has only finite number of
fixed points, it is clear that $A(\phi,J_0) > 0$ and so we have $
A(\phi) > 0$. Note that when $(M,\omega)$ is weakly exact and so
$A_{S}(\phi,J_0;J') = \infty$, $A(\phi,J_0)$ is reduced to
$$
A(\phi,J_0) = \sup_{J\in j_{(\phi,J_0)}}\{ A_D(\phi,J_0;J') \}.
$$

In [Oh5] we have proved
$$
A(\phi) \leq \gamma(\phi) \leq \|\phi\|_{mid}, \tag 1.17
$$
postponing to the present paper the proof of some existence
theorem of the following equation
$$
\cases
\dudtau + J_t\Big(\dudt - X_H(u)\Big) = 0 \\
u(-\infty) = [z^-,w^-], \, u(\infty) = [z^+,w^+] \\
w_- \# u \sim w_+, \quad u(0,0)=q \in B(u)
\endcases
\tag 1.18
$$
for the particular family
$$
J_t = (\phi_H^t)_*J'_t, \quad J' \in j_{(\phi,J_0)}. \tag 1.19
$$
{\it Note that $J_t$ in (1.19) is $t$-periodic by definition of
$j_{(\phi,J_0)}$}, i.e. $J_1 = J_0$ and hence (1.18) is
well-posed.

Let $k: M \to \R$ be a Morse function and $\e > 0$ be any small
positive number. We denote by $\alpha, \, \beta \in CF_n(\e k)$
Floer cycles of $\e k$ that realize the fundamental cycle $1^\flat
= [M]$ and $h_\HH(\alpha), \, h_\HH(\beta)$ be the cycles in
$CF_n(H)$ of $H$ transferred via the linear homotopy
$$
\HH: s\in [0,1] \mapsto (1-s)\e k + s H.
$$
What we needed in the proofs of the optimal Energy-Capacity
inquality and (1.17) in [Oh5] is the following
existence theorem of {\it non-stationary} solutions of (1.18) {\it
with the upper bound on the action}. We refer to [Oh4,5] for more
explanations on the definitions and notations of various terms
undefined in this theorem. The proof of this theorem (Theorem 7.1)
will be finished in section 7.

\proclaim{Theorem II [Theorem 3.11, Oh5]} Let $H$ and $J_0$ be as
before. And let $q \in \text{Int }B(u)$ and $\delta
>0$ be given. Then for any $J' \in j_{(\phi,J_0)}$, there exist
some generators $[z,w] \in h_{\HH}(\alpha)$ and $[z', w'] \in
h_{\HH}(\beta)$ with
$$
\aligned \AA_H([z,w]) & \leq \rho(H;1) + \frac{\delta}{2} \\
\AA_{\widetilde H}([z',w']) & \leq  \rho(\widetilde H;1) +
\frac{\delta}{2}
\endaligned
\tag 1.20
$$
such that the following alternative holds: \roster
\item The equation
$$
\cases \dudtau + J  \Big(\dudt - X_H(u)\Big) = 0 \\
u(\infty)=[\widetilde z',\widetilde w'], \quad u(-\infty) = [z,w]
\endcases
$$
has a cusp-solution
$$
u_1 \# u_2 \# \cdots \cdots \# u_N
$$
which is a connected union of Floer trajectories, possibly with
a finte number of sphere bubbles, for $H$ that
satisfies the conditions
$$
u_N(\infty) = [\widetilde z',\widetilde w'], \, u_1(-\infty) =
[z,w], \quad u_j(0,0) = q \in B(u).
$$
for some $1 \leq j \leq N$,
\item or there is some $J'_t$-holomorphic sphere $v:S^2 \to M$
 for some $t\in [0,1]$ that passes through the point $q$.
\endroster
This in particular implies
$$
0 < A(\phi,J_0)\leq A(\phi)\leq \gamma(\phi) < \infty
$$
for any $\phi$ and $J_0$.
\endproclaim

We would like to emphasize that without the upper estimate (1.20)
the existence of such pair $[z,w]$ and $[z',w']$ would have been a
result much easier to prove from the nontriviality of the quantum
product
$$
1 * 1 = 1
$$
and its interpretation in terms of the pants product in the Floer
complex. However obtaining the energy bound (1.20) requires
careful control of the levels of the Floer cycles representing the
Floer class $1^\flat$ under the pants product. We have carried out
this analysis exploiting the general properties of spectral
invariants $\rho(\cdot;1)$ and the energy identity that we derive
for the deformed mapping cylinder in section 3.

Using our techniques proving this existence result, we can in fact
prove a more general existence theorem (see Theorem 5.1) on
certain piecewise smooth finite energy solutions, i.e. solutions
satisfying
$$
\int \int \Big(\Big| \dudtau\Big|^2 + \Big|\dudt -
X_{(H,F)}\Big|^2\Big) dt\,d\tau  < \infty,
$$
of the following perturbed Cauchy-Riemann equation with {\it
discontinuous} Hamiltonian perturbation term
$$
\cases \dudtau + J \Big(\dudt - X_{(H,F)}(u) \Big)= 0 \\
u(\infty)=z^+ \in \text{Per}(F), \quad u(-\infty) = z^- \in
\text{Per}(H)
\endcases
\tag 1.21
$$
where $X_{(H,F)}$ is the discontinuous family of vector field
$$
X_{(H,F)}(\tau,t,x) = \cases X_H(t,x) \quad \text{for } \tau < 0 \\
X_F(t,x) \quad \text{for }\tau > 0
\endcases
\tag 1.22
$$
and $J$ is the discontinuous family
$$
J(\tau,t,x) =\cases (\phi_H^t)_*J'_t \quad\text{for } \, \tau < 0 \\
(\phi_F^t)_*J'_t \quad\text{for } \, \tau > 0
\endcases
\tag 1.23
$$
for $J' \in j_{(\phi,J_0)}$. Note that if we define
$$
v(\tau,t):= \cases (\phi_H^t)^{-1}(u(\tau,t)) \quad\text{for }\,
\tau < 0 \\
(\phi_F^t)^{-1}(u(\tau,t)) \quad\text{for }\, \tau > 0
\endcases
$$
then $v$ satisfies (1.9) away from $\tau = 0$ and has finite
energy i.e., lie in $W^{1,2}$. A priori, $v$ may be discontinuous
but since $J'_t$ is smooth everywhere, $W^{1,2}$ condition implies
that $v$ must be smooth across $\tau=0$ and becomes the classical
solution of (1.9), if it is continuous.

We refer to \S 4 and \S 5 for more detailed study of (1.21),
especially Theorem 5.1. Theorem 5.1 is the main existence result
of the present paper, whose statement however is somewhat long and
awkward to state in this introduction. We refer readers to section
5 for the precise statement. We would like to point out that
Kasturirangan and the author [KO] already demonstrated the
necessity of studying the Cauchy-Riemann equation with {\it
non-smooth} Lagrangian boundary condition. The appearance of the
perturbed Cauchy-Riemann equation (1.20) with discontinuous
Hamiltonian perturbation term is somehow reminiscent of the story
from [KO]. Both appear in our study of the {\it chain level}
operators in the Floer theory when we go to some limiting cases
where the usual Floer's perturbed Cauchy-Riemann equation is not
analytically well-posed but allows smooth approximations which
enable us to take the limit of the homology. The existence of this
limiting homology (`Fary functor') is the source of this existence
theorem of discontinuous $W^{1,2}$ solution of the above equation.

One motivation of ours to study (1.21) is that this equation is
naturally related to the {\it small Hofer pseudo-norm}: the small
Hofer pseudo-norm, denoted by $\|\cdot\|_{sm}$ is defined by
$$
\|\phi\|_{sm}:= \rho^+(\phi) + \rho^-(\phi)
$$
where
$$
\align
\rho^+(\phi) & = \inf_{H \mapsto \phi} \int_0^1 \max H_t\, dt
= \inf_{\pi(h) = \phi}E^+(h)\\
\rho^-(\phi) & = \inf_{F \mapsto \phi} \int_0^1 -\min F_t\, dt
= \inf_{\pi(h) = \phi}E^-(h)
\endalign
$$
The question whether this pseudo-norm is nondegenerate was first posed by
Polterovich [Po] and proved by McDuff [Mc] for the case of
$\C P^n$ or for the weakly exact case.
We refer to [Mc], [Po] for some more background materials on various
related norms of Hamiltonian diffeomorphisms.

\proclaim{Theorem III} Let $X_{(H,F)}$ and $J$ as in (1.21) and
(1.22) respectively. Suppose $u: \R \times S^1 \to M$ is a
piece-wise smooth map with finite energy that satisfies
$$
\cases \dudtau + J \Big(\dudt - X_{(H,F)}(u)\Big) = 0 \\
u(\infty)=[z^+,w^+] \in \text{Crit}\AA_F, \quad u(-\infty) =
[z^-,w^-] \in \text{Crit}\AA_H
\endcases
\tag 1.24
$$
and
$$
\int (w^-)^*\omega + \int (u_-)^*\omega + \int (u_+)^*\omega -
\int (w^+)^*\omega = 0 \tag 1.25
$$
where $u_\pm$ are the parts of $u$ on $\tau <0$ and $\tau > 0$
respectively. Then we have the following inequality
$$
\aligned \int \Big| \dudtau \Big|^2_{J_t}
& \leq  -\AA_F([z^+,w^+]) + \AA_H([z^-,w^-]) \\
& \quad + \int_0^1 \max H\, dt + \int_0^1 -\min F_t\, dt.
\endaligned
\tag 1.26
$$
\endproclaim
This is the key inequality which relates all three basic
quantities in the Floer theory (or more generally in the
symplectic topology), the {\it energy} of the pseudo-holomorphic
curves, the {\it actions} of periodic orbits and the {\it Hofer
type quantities} of the corresponding Hamiltonians. We hope to pursue
applications of Theorem 5.1 and Theorem III in the future.

Finally we would like to note that the kind of singular
degeneration problem that we consider in section 5 and 6 will be
the first step towards understanding the complete picture of
adiabatic degeneration of Floer's perturbed Cauchy-Riemann
equation which will be studied jointly with Fukaya [FOh2].

We would like to thank the Institute for Advanced Study in
Princeton for the excellent environment and hospitality during our
participation of the year 2001-2002 program ``Symplectic Geometry
and Holomorphic Curves''. Much of this research has been carried
out while we are visiting the Korea Institute for Advanced Study
in Seoul. We thank KIAS for providing excellent atmosphere of
research.

\head{\bf \S 2. Perturbed Cauchy-Riemann equations re-visited}
\endhead

In this section, we re-examine Floer's perturbed Cauchy-Riemann
equation that defines the chain map in the Floer homology theory
and show how the Cauchy-Riemann equation with jumping
discontinuity in its coefficients naturally appears in the chain
level Floer theory.

The story goes in the following way. Let $H, \, F$  be two
$t$-periodic Hamiltonians and $J^0, \, J^1$ with $J^i=\{J^i_t\}_{0
\leq t \leq 1} , \, i = 0,\, 1$ be a $t$-periodic family of almost
complex structures. Without loss of any generality, we will always
assume
$$
H_t \equiv 0, \quad J_t \equiv J_0 \quad\text{near }\, t = 0\equiv
1 \tag 2.1
$$
When we are given a homotopy $\HH$ with $\HH = \{H^s\}_{0\leq s
\leq 1}$ with $H^0 = H, \quad H^1=F$ and $j=\{J^s\}_{0\leq s\leq
1}$ connecting $J^0$ and $J^1$, the chain homomorphism
$$
h^\rho_{(\HH,j)}: CF_*(H,J^0) \to CF_*(F,J^1)
$$
is defined by the non-autonomous equation
$$
\cases \frac{\part u}{\part \tau} +
J^{\rho_1(\tau)}\Big(\frac{\part u}{\part t}
- X_{H^{\rho_2(\tau)}}(u)\Big) = 0\\
\lim_{\tau \to -\infty}u(\tau) = z^-, \,  \lim_{\tau \to
\infty}u(\tau) = z^+
\endcases
\tag 2.2
$$
with finite energy condition
$$
\int \Big| \dudtau \Big|_{J^{\rho_1(\tau)}_t}^2 < \infty. \tag 2.3
$$
Here $\rho_i$ are the cut-off functions of the type $\rho :\R \to
[0,1]$,
$$
\align
\rho(\tau) & = \cases 0 \, \quad \text {for $\tau \leq R_1$}\\
                    1 \, \quad \text {for $\tau \geq R_2$}
                    \endcases \\
\rho^\prime(\tau) & \geq 0
\endalign
$$
for given arbitrary pair $R_1 < R_2$ in $\R$. For the simplicity
of exposition, we will consider mostly the case
$\rho_1=\rho_2=\rho$. (2.2) and (2.3) uniquely determine the
asymptotic condition
$$
\lim_{\tau \to -\infty}u(\tau) = z^- \in \text{Per}(H), \, \,
\lim_{\tau \to \infty}u(\tau) = z^- \in \text{Per}(F) \tag 2.4
$$
Since we consider (2.2) as an equation on the $\Gamma$-covering
space $\widetilde \Omega_0(M)$, we also put the topological
condition
$$
w^- \# u \sim w^+ \tag 2.5
$$
where $w^\pm$ are discs bounding $z^\pm$ respectively. Here
$\Gamma$ is the covering group
$$
\Gamma = {\pi_2(M) \over \ker (\omega|_{\pi_2(M)}) \cap \ker
(c_1|_{\pi_2(M)})}.
$$
Then we lift the asymptotic condition (2.4) to the covering space
as
$$
\lim_{\tau \to -\infty}u(\tau) = [z^-, w^-]  \in \text{Crit}\AA_H,
\, \,  \lim_{\tau \to \infty}u(\tau) = [z^+,w^+] \in
\text{Crit}\AA_F. \tag 2.6
$$
The following identity is the fundamental identity which has been
used frequently in our previous works [Oh2,3,4]. We leave its
proof to readers.

\proclaim{Lemma 2.1} For any solution of (2.2) with (2.5), we have
$$
\align \AA_F([z^+,w^+]) & - \AA_H([z^-,w^-]) \\
& = - \int \Big|\dudtau \Big|_{J^{\rho_1(\tau)}}^2 -
\int_{-\infty}^\infty \rho'_2(\tau)(F(t,u(\tau,t)) -
H(t,u(\tau,t))) \, dt\,
d\tau \tag 2.7 \\
& \leq - \int \Big|\dudtau \Big|_{J^{\rho_1(\tau)}}^2 + \int_0^1
-\min F_t\, dt + \int_0^1 \max H_t\, dt.
\endalign
$$
In particular, we have
$$
 \int \Big|\dudtau
\Big|_{J^{\rho_1(\tau)}}^2 \leq \AA_H([z^-,w^-]) -\AA_F([z^+,w^+])
+ \int_0^1 -\min F_t\, dt + \int_0^1 \max H_t\, dt. \tag 2.8
$$
\endproclaim
We recall that when $H=F$ and $\rho \equiv 1$, we have the
improved identity
$$
\int \Big|\dudtau \Big|_J^2 = \AA_H([z^-,w^-]) -\AA_H([z^+,w^+]).
$$
Note that the upper bound (2.8) is independent of the cut-off
functions $\rho$ and also of $R_1, \, R_2$. One important fact
which we are going to exploit is that the chain map itself
$$
h^\rho_{(\HH,j)}: CF_*(H,J^0) \to CF_*(F,J^1)
$$
depends on the choice of $\rho$ although two chain maps
$h^\rho_{(\HH, j)}$ for different $\rho$'s are chain homotopic to
each other and hence they induce the same homomorphism
$$
h_{HF}: HF_*(H) \to HF_*(F)
$$
in homology.

We now would like to study the equation (2.2) with the  condition
(2.5) as $R_1 \to 0_-$ and $R_2 \to 0_+$, i.e., when the cut-off
function $\rho$ converges to the {\it discontinuous} Heaviside
function
$$
\rho =\cases 0 \quad \text{for } \, \tau < 0 \\
1 \quad \text{for } \, \tau > 0.
\endcases
$$
Note that the corresponding limit of the equation (2.2) is
$$
\cases \frac{\part u}{\part \tau} +
J\Big(\frac{\part u}{\part t} - X_{(H,F)}(u)\Big) = 0\\
\lim_{\tau \to -\infty}u(\tau) = z^-,  \, \lim_{\tau \to
\infty}u(\tau) = z^+
\endcases
\tag 2.10
$$
where $X_{(H,F)}$ are the {\it discontinuous} family of
Hamiltonian vector fields satisfying (1.21)  and
$$
J(\tau,t,x) =\cases J^0(t,x) \quad\text{for } \, \tau < 0 \\
J^1(t,x) \quad\text{for } \, \tau > 0 \endcases.
$$
We will be
particularly interested in the family
$$
J^0_t = (\phi_H^t)^*J'_t, \quad J^1_t = (\phi_F^t)^*J'_t \tag 2.11
$$
for $J' \in j_{(\phi,J_0)}$. We note that {\it these are
$t$-periodic} due to the definition of $j_{(\phi,J_0)}$ and so
(2.10) is a well-defined equation over $\R \times S^1$ in the
classical sense, except at $\tau = 0$. Due to the condition (2.1),
the equation has smooth coefficients and the Hamiltonian
perturbation term is smooth in a neighborhood of $\R \times \{0\}
\subset \R \times S^1$ including $(0,0)$. Therefore any $W^{1,2}$
solution of (2.10) will be smooth near $(0,0)$. More importantly,
under the choice (2.11), (2.10) is transformed into
$$
\cases {\part v\over \part \tau} + J'_t{\part v \over \part t} = 0 \\
\phi(v(\tau,1)) = v(\tau,0)
\endcases
\tag 2.12
$$
by the map
$$
v(\tau,t) = \cases (\phi_H^t)^{-1}(u(\tau,t)) \quad \text{for } \,
\tau < 0 \\
(\phi_F^t)^{-1}(u(\tau,t)) \quad \text{for } \, \tau > 0.
\endcases
\tag 2.13
$$
If $u$ is in $W^{1,2}$, so is $v$. Since (2.12) has smooth
coefficients, such $W^{1,2}$ solution $v$ will be indeed smooth if
$v$ is continuous across $\tau = 0$.

The equation (2.10) now has bounded discontinuous coefficients and
zero-order perturbations. Therefore we regard (2.10) as the first
order quasi-linear PDE for the {\it distributional} maps. For this
purpose, we will always embed $M$ with a fixed metric $g$ into
$\R^N$ by the Nash isometric embedding theorem and regard $u$ as a
vector valued distributional maps defined on $\R \times S^1$. As
usual in the geometric PDE, we define a {\it finite energy} map
into $M$ by a $\R^N$-valued distribution $u$ with
$$
u(z) \in M \subset \R^N \quad \text{almost everywhere}
$$
and
$$
\int \Big| \dudtau\Big|^2 + \Big|\dudt - X_{(H,F)}\Big|^2
dt\,d\tau < \infty \tag 2.14
$$
which does not depend on the choice of metric $g$. Note that since
$M$ is compact, (2.14) in particular implies that $u \in L^\infty$
and so in $L^p$ for every $p > 0$.

The following is the key property that will be important later in
the nondegeneracy proof.

\proclaim{Theorem 2.2} Let $u: \R \times S^1 \to M$ be a solution
lying with finite energy which satisfies (2.10) and (2.6). Suppose
that it satisfies in addition
$$
\int w_-^*\omega + \int u_-^*\omega + \int u_+^*\omega- \int
w_+^*\omega = 0. \tag 2.15
$$
Then we have the following  identity
$$
\aligned \int \Big| \dudtau \Big|^2_{J_t} & \leq -
\AA_F([z^+,w^+]) + \AA_H([z^-,w^-]) \\
& \quad + \int_0^1 \max H_t \, dt + \int_0^1 -\min F_t\, dt.
\endaligned
\tag 2.16
$$
\endproclaim
\demo{Proof} We first note that since $u$ has finite energy, both
$\int u_+^*\omega$ and $\int u_-^*\omega$ are finite and so
well-defined (see [Oh2] for the proof of this fact). It also
implies
$$
\align \int \Big| \dudtau \Big|_{J}^2 & = \lim_{\e \to 0}
\int_{(\R \setminus [-\e,\e]) \times S^1} \Big| \dudtau \Big|_{J}^2 \\
& = \lim_{\e \to 0} \int_{-\infty}^{-\e}\int_0^1 \Big| \dudtau
\Big|_{J^0_t}^2 \, dt\, d\tau + \lim_{\e \to 0}
\int_\e^\infty\int_0^1 \Big| \dudtau \Big|_{J^1_t}^2\, dt\, d\tau.
\endalign
$$
On the other hand, we have
$$
\align \int_{-\infty}^{-\e}\int_0^1 \Big| \dudtau \Big|_{J^0_t}^2
\, dt\, d\tau
& = \AA_H(u(-\e)) - \AA_H(u(-\infty)) \\
\int_{-\infty}^{-\e}\int_0^1 \Big| \dudtau \Big|_{J^1_t}^2 \, dt\,
d\tau & = \AA_F(u(\infty)) - \AA_F(u(+\e))
\endalign
$$
and hence
$$
\int \Big|\dudtau\Big|^2_{J} = -\AA_F(u(\infty)) +
\AA_H(u(-\infty)) + \lim_{\e \to 0}(-\AA_H(u(-\e)) + \AA_F(u(\e)).
\tag 2.17
$$
It follows from (2.15) that the term of the integral of $\omega$
in the action functional vanishes and so
$$
\align \lim_{\e \to 0} ( -\AA_H(u(-\e)) + \AA_F(u(\e)) & =
\lim_{\e \to 0}
\Big(\int_0^1 -F(u(\e,t))\, dt + \int H(u(-\e,t))\, dt\Big) \\
& \leq \int_0^1 \max H_t \, dt + \int_0^1 -\min F_t\, dt. \tag
2.18
\endalign
$$
Combining (2.17) and (2.18), we have finished the proof.
\qed\enddemo

The essential matter now is to prove some existent result of the
equation (2.2) with suitable asymptotic condition and  (2.15). One
way of obtaining a solution is by taking a limit of a sequence of
solutions of (2.2) as $\rho$ converges to the step function
(2.10). But this limit could be {\it stationary}, i.e,
$\tau$-independent which is not suitable for our purpose because
the solution we need is the one that is non-trivial as in [Oh5].
For this purpose, we need to generalize the above story to the
maps from the arbitrary compact Riemann surface of genus zero with
$k$ punctures. We will focus only on $k=3$  in this paper because
this is the case that is directly relevant to the proof of Theorem
I.

\head{\bf \S3. Energy identity on the deformed mapping cylinder}
\endhead

In this section, we recast Floer's equation for the chain map in
the point of the deformed mapping cylinder, and calculate the
vertical energy  of the pseudo-holomorphic sections for a suitably
chosen almost complex structure $\widetilde J$ in terms of the
variation of the Hamiltonians and the curvature of the natural
connection associated to the Hamiltonians. This will form the
basis of the energy calculation of  pseudo-holomorphic sections of
the Hamiltonian fibration over the Riemann surface $\Sigma$ of
genus zero with arbitrary number of punctures, when we realize the
conformal structure of $\Sigma$ by the minimal area metric on
$\Sigma$ [Z]. Our calculation of the  energy is a crucial
ingredient in the optimal Energy-Capacity inequality and
nondegeneracy proof of $\gamma$-norm,  and also in
our construction of $W^{1,2}$-solution mentioned in Theorem II (or
Theorem 5.1). This calculation is very much in the spirit of geometric
calculations in the differential geometry.

Consider the two parameter family of Hamiltonian diffeomorphisms
$$
\phi: (s,t) \in [0,1] \times [0,1] \to \HH am(M,\omega)
$$
with
$$
\phi(s, 0) = id
$$
for all $s \in [0,1]$. We also denote $\phi_s^t = \phi(s,t)$ and
$$
\phi_s: = \phi(s,1).
$$
We will assume that
$$
\phi(s,t)  \equiv \cases id \quad \text{for $t$ near 0}\\
\phi_s \quad \text{for $t$ near 1}
\endcases
$$
We define the vector fields
$$
\aligned
X(s,t,x) & = {\part \phi_s^t \over \part t} \circ (\phi_s^t)^{-1} \\
Y(s,t,x) & = {\part \phi_s^t \over \part s} \circ (\phi_s^t)^{-1}
\endaligned
\tag 3.1
$$
and denote by
$$
\align
H & : [0,1] \times [0,1] \times M \to \R\\
F & : [0,1] \times [0,1] \times M \to \R
\endalign
$$
the {\it normalized} Hamiltonians generating $X$ and $Y$
respectively. We will be particularly interested in the family
generated by the {\it linear} homotopy
$$
H(s,t,x) = s H(t,x).
$$
Now we consider the cylinder $\Sigma = \R \times S^1$ or the
semi-cylinders
$$
\align
\Sigma_+ & = [0,\infty) \times S^1 \\
\Sigma_- & = (-\infty,0] \times S^1
\endalign
$$
equipped with the standard flat metric on them. We will focus on
$\Sigma=\Sigma_-$ but the parallel story goes for $\Sigma$ or
$\Sigma_+$ with some obvious changes. For each given $0 < R_1 <
R_2$, we consider the cut-off functions
$$
\rho= \rho^-_{R_1,R_2} = \cases 0 \quad \text{for } \, -\infty
< \tau \leq - R_2 \\
1 \quad \text{for } \, -R_1 \leq \tau \leq 0
\endcases
$$
with $\rho' \leq 0$ using the function $\rho$. We reparameterize
the $s$-variable in the family $\phi$ and denote
$$
\phi_\rho(\tau,t) = \phi(\rho(\tau),t).
$$
We denote by $X_\rho$ and $Y_\rho$ the corresponding
reparameterized family of vector fields. If we denote by $H_\rho$
and $F_\rho$ the corresponding Hamiltonians for the
reparameterized family $\phi_\rho$, the following is easy to check
$$
\aligned
H_\rho(\tau,t,x) & = \rho(\tau) H(t,x) \\
F_\rho(\tau,t,x) & = \rho'(\tau) F(\rho(\tau),t,x).
\endaligned
\tag 3.2
$$
Considering the path
$$
f: \tau \mapsto f_\tau=\phi^1_{\rho(\tau)}\, ; \quad \R_- \to \HH
am(M,\omega)
$$
we define the {\it deformed mapping cylinder} by
$$
E = E_f := \R_- \times \R \times M / (\tau,t,x) \sim
(\tau,t-1,f_\tau(x))
$$
which defines a symplectic fibration over the semi-cylinder
$\Sigma= \R_+ \times S^1$. We consider the closed two form
$$
\omega_E:= \omega + d(H_\rho dt) = \omega + d(\rho H_t dt) \tag
3.3
$$
on $\R_- \times \R \times M$ which projects down to  a closed two
form on $E$. We denote the push down of the form again by
$\omega_{E}$. This form is nondegenerate in the fiber of $E \to \R
\times S^1$ and restricts to $\omega$ and so defines a canonical
connection $\nabla$ [GLS]  whose coupling form is exactly
$\omega_{E}$. We need to describe this connection more explicitly
to carry out precise calculation of the energy.

We trivialize the mapping cylinder $E\to \R_-\times S^1 \times M$
by
$$
[\tau,t,v] \in E \mapsto (\tau,t,u(\tau,t)) \tag 3.3
$$
where
$$
u (\tau,t) = (\phi^t_{\rho(\tau)})^{-1}(v). \tag 3.4
$$
Note that (3.4) is $t$-periodic because of the defining
equivalence relation of the deformed mapping cylinder $E_f$ and so
(3.4) provides a well-defined trivialization of $E_f$.

Under this trivialization, a straightforward calculation proves
the formula for the horizontal lifts of ${\part \over \part \tau}$
and ${\part \over \part t}$ for the connection $\nabla$
$$
\aligned
(Dv)^h({\part \over \part \tau}) & = {\part \over \part \tau}\\
(Dv)^h({\part \over \part t} ) & = {\part \over \part t} + X_\rho
\endaligned
\tag 3.5
$$
for the arbitrary sections $v$. If we denote by
$$
\Pi_{\tau,t,x}: T_{(\tau,t,x)}E \to (TE)^v_{(\tau,t,x)} \cong T_xM
$$
the associated vertical projection, it has the formula
$$
\Pi_{\tau,t,x}(\alpha,\beta,\xi) = \xi - \beta X_\rho \tag 3.6
$$
in the trivialization (3.3). We provide a symplectic form on $E$
by
$$
\Omega_\lambda= \Omega_{E,\lambda} = \omega_E + \lambda
\omega_\Sigma
$$
for a sufficiently large $\lambda >0$, where $\omega_\Sigma$ is an
area form on the base $\R_+ \times S^1$. We normalize it so that
$$
\int \omega_\Sigma = 1. \tag 3.7
$$

Now let $J^s_t$ for $(s,t) \in [0,1]^2$ be two parameter family of
$\omega$-compatible almost complex structures on $M$ such that $J$
is constant near $s = 0, \, 1$ and $J^s_0 = J^s_1$. We also assume
that $J^s_t$ is constant near $t = 0 \equiv 1$. Using this, we
define the almost complex structure $\widetilde J$ on $E$ in the
trivialization (3.3) by the formula
$$
\widetilde J(\tau,t,x)(\alpha,\beta,\xi) = (-\beta,\alpha,
(\phi_{\rho(\tau)}^t)^*J_t^{\rho(\tau)} (\xi - \beta X_\rho) +
\alpha X_\rho). \tag 3.8
$$
One can easily check that $\widetilde J$ is well-defined on $E$
and tame to $\Omega_{E}$ for a sufficiently large $\lambda$, but
it is {\it not compatible} in the usual sense in that the bilinear
form
$$
\Omega(\cdot, \widetilde J \cdot)
$$
is not symmetric. However we can symmetrize this and define the
associate metric $g_{\widetilde J}$ by
$$
\langle V,W \rangle =  g_{\widetilde J}(V,W): = {1 \over
2}(\Omega(V, \widetilde J W) + \Omega(W, \widetilde J V)). \tag
3.9
$$
We call (3.9) the metric associated to $\widetilde J$ and denote
$$
|V|^2 = |V|^2_{\widetilde J} = g_{\widetilde J}(V,V).
$$
With respect to this metric, we still have the following basic
identity whose proof we omit.

\proclaim{Lemma 3.1} Let $v: \R \times S^1 \to E$ be any
$\widetilde J$-holomorphic map $v$. Then we have
$$
{1 \over 2} \int |Dv|^2_{\widetilde J} = \int v^*\Omega_\lambda. \tag 3.11
$$
\endproclaim

Next note that for any $\widetilde J$-holomorphic section $v:
\Sigma \to E_{f}$, if we identify it with a map
$$
v: \R_+ \times \R \to M
$$
satisfying
$$
f_{\tau}(v(\tau,1)) = v(\tau,0),
$$
the map $u$ defined by
$$
u(\tau,t) = (\phi_{\rho(\tau)}^t)^{-1}(v(\tau,t))
$$
satisfies
$$
u(\tau,0) = u(\tau,1)
$$
and is smooth near $t = 0\equiv 1$ because of the condition (2.1)
for the Hamiltonian $H$, and hence defines a well-defined map $u:
\R_+ \times S^1 \to M$. It also satisfies
$$
\dudtau + J^{\rho(\tau)}_t\Big(\dudt -
X_{H_t^{\rho(\tau)}}(u)\Big) = 0. \tag 3.12
$$
With this preparation, we now compute the energy density $|Dv|^2$
for $\widetilde J$-holomorphic section $v : \Sigma \to E_{f}$.
Decomposing $Dv$ into its vertical and horizontal components
$$
Dv = (Dv)^v + (Dv)^h
$$
we have
$$
|Dv|^2 = |(Dv)^v|^2 + |(Dv)^h|^2 + 2 \langle (Dv)^v, (Dv)^h
\rangle. \tag 3.13
$$
If we denote  the curvature of the connection $\nabla$ by
$$
K(v) d\tau \wedge dt
$$
then it is straightforward to compute
$$
{1\over 2} \int |(Dv)^h|^2 = \int K(v) d\tau \,dt + \lambda \tag
3.14
$$
by integrating the identity
$$
\align \sum_{i=1}^2|(Dv)^h(e_i)|^2 & =
\sum_{i=1}^2\Omega_{E,\lambda}((Dv)^h(e_i),
\widetilde J (Dv)^h(e_i)) \\
& = \sum_{i=1}^2 (\omega_E + \lambda \omega_{\Sigma})((Dv)^h(e_i),
\widetilde J (Dv)^h(e_i)) \\
& = 2 (\omega_E((Dv)^h(e_1), (Dv)^h(e_2)) + \lambda
\omega_\Sigma(e_1,e_2))
\endalign
$$
for an orthonormal frame $\{e_1, e_2\}$ of $T\Sigma$ , e.g., for
the frame
$$
\Big\{{\part \over \part \tau}, {\part \over \part t} \Big\}
$$
and then applying the curvature identity
$$
d(\omega(e_1^\#,e_2^\#)) = -\iota_{[e_1^\#,e_2^\#]}  \omega
$$
where $e_i^\#$ is the horizontal lift of $e_i$ (see [(1.12), GLS]
but with caution on the sign convention). Applying (3.5) and
(3.8), we derive
$$
\aligned \widetilde J (Dv)^h({\part \over \part \tau}) & = {\part
\over \part t}
+ X_\rho \\
\widetilde J (Dv)^h({\part \over \part t}) & = -{\part \over \part
\tau}.
\endaligned
\tag 3.15
$$
(3.5) immediately gives rise to
$$
K(v)(\tau,t) = \omega_E\Big((Dv)^h({\part \over \part \tau}),
(Dv)^h({\part \over \part t})\Big) = \rho'(\tau) H(t, u(\tau,t))
\tag 3.16
$$
which in particular proves that $K$ is compactly supported and so
the curvature integral in (3.14) in finite. Combining (3.14) and
(3.16), we have derived
$$
{1\over 2} \int |(Dv)^h|^2 = \lambda + \int \rho'(\tau) H(t,
u(\tau,t)) \tag 3.17
$$
Next we will compute $\langle (Dv)^v, (Dv)^h\rangle$
$$
\langle (Dv)^v, (Dv)^h\rangle = \langle (Dv)^v({\part \over
\part \tau}), (Dv)^h({\part \over \part \tau})\rangle + \langle
(Dv)^v({\part \over \part t}), (Dv)^h({\part \over \part t})
\rangle.
$$
We separately compute the two terms. First from the definition
(3.9) of the associated metric, we have
$$
\align \langle (Dv)^v({\part
\over \part \tau}), & (Dv)^h({\part \over \part \tau})\rangle \\
& = {1 \over 2}\Big(\Omega((Dv)^v({\part \over \part \tau}),
 \widetilde J(Dv)^h({\part \over \part \tau}))  +
\Omega((Dv)^h({\part \over \part \tau}), \widetilde J
(Dv)^v({\part \over \part \tau}))\Big).
\endalign
$$
Since $\widetilde J$ preserves the vertical space, it is easy to
see that the second term vanishes. For the first term, we derive
from (3.6) and (3.15)
$$
2 \Omega((Dv)^v({\part \over \part \tau}), \widetilde J
(Dv)^h({\part \over
\part \tau})) = 2 \omega_E(\dudtau, {\part \over
\part t} + X_\rho) = 0
$$
and hence
$$
2 \langle (Dv)^v({\part \over \part \tau}), (Dv)^h({\part \over
\part \tau})\rangle = 0
\tag 3.18
$$
On the other hand, a straightforward calculation using (3.6) and
(3.15) also show
$$
\align 2 \langle (Dv)^v({\part \over \part t}), &(Dv)^h({\part
\over \part t})\rangle \\
& = \Big(\omega + d(\rho H_t dt)\Big) (\dudt - X_\rho, -{\part
\over \part \tau}) = 0. \tag 3.19
\endalign
$$
Adding (3.18) and (3.19), we derive the identity
$$
2 \langle (Dv)^v, (Dv)^h \rangle  = 0. \tag 3.20
$$
We summarize the above calculations into

\proclaim{Lemma 3.3} Let $\widetilde J$ be the above almost
complex structure on $E$ and $v: \Sigma \to E$ be a $\widetilde
J$-holomorphic section. Then we have
$$
\align {1\over 2} \int |(Dv)^h|^2 & = \lambda + \int \rho'(\tau)
H(t, u(\tau,t)) \\
\langle (Dv)^v, (Dv)^h \rangle &  = 0
\endalign
$$
\endproclaim

Applying Lemma 3.3, we derive the following formula for the
vertical energy

\proclaim{Proposition 3.4} Let $\rho = \rho_{R_1,R_2}$ be the
cut-off functions as above. Then we have
$$
\align {1 \over 2} \int |(Dv)^v|^2 & = \int v^*\Omega_\lambda - {1
\over 2} \int |(Dv)^h|^2
- \int \langle (Dv)^v, (Dv)^h \rangle \tag 3.21\\
& = \int v^*\omega_{E} -  \int \rho'(\tau) H(t, u(\tau,t))\, dt\,
d\tau \tag 3.22 \\
\endalign
$$
for any $\rho$ and $\widetilde J$-holomorphic section $v:
\Sigma_-=\R_- \times S^1 \to E$ of the deformed mapping cylinder.
Similar formula holds on $\Sigma$ or $\Sigma_+ = \R_+\times S^1$.
In particular, we have the inequality
$$
{1 \over 2} \int |(Dv)^v|^2 \leq \int v^*\omega_{E} + \int_0^1
\max H_t\, dt \tag 3.23
$$
for any finite energy $\widetilde J$-holomorphic sections $v$ for
any $0 < R_1 <R_2 < \infty$.
\endproclaim
\demo{Proof} (3.22) immediately follows by substituting (3.16) and
(3.20) into (3.21). (3.23) is a consequence of integration over
$\tau$ noting that $\rho' \leq 0$ on $\Sigma_- = (-\infty,0]
\times S^1$. \qed\enddemo

\head{\bf \S4. Pants product and the Hamiltonian fibration}
\endhead

Let $h, f, g \in \widetilde{\HH am}(M,\omega)$ be an arbitrary
nondegenerate triple.  The product structure
$$
HF_*(h) \otimes HF_*(f) \to
HF_*(g). \tag 4.1
$$
is defined by considering the ``pants product'' in the chain complex
$$
CF_*(H) \otimes CF_*(F) \to CF_*(G)
\tag 4.2
$$
where $[H] = h, \, [F] =f, \, [G] = g$. Here one can already see
that the pants product in the chain level depend on the choice of
the Hamiltonians representing the classes  $h,\, f,\, g \in
\widetilde{\HH am}(M,\omega)$, not just on $h,\, f, \, g$ let
alone their projections $\phi=\pi(h), \, \psi=\pi(f), \, \eta =
\pi(g)$. The definition of the pants product (4.2)  depend on many
additional data, while the induced product in homology (4.1) is
independent of the choice of such data. The best way of describing
the pants product is using the Hamiltonian fibration with
connection of fixed boundary monodromy as described by Entov [En].
We briefly recall Entov's construction here with few notational
changes and different convention on the grading, but fully
exploiting the {\it minimal area metric representation} [Z] of the
conformal structure on $\Sigma$. As in [Oh5], our grading will be
provided by
$$
k = \mu_H 
$$
where $\mu_H([z,w])$ is the Conley-Zehnder index of the critical point
$[z,w] \in \AA_H$. This grading  convention respects
the grading under the pants product
$$
CF_k(H) \otimes CF_\ell(F) \to CF_{k + \ell-n}(G).
$$

Let $\Sigma$ be the compact Riemann surface of genus 0 with three
punctures with the minimal area metric [Z] which we now describe.
We conformally identify $\Sigma$ with three half cylinders  which
we denote by $\Sigma_1$, $\Sigma_2$ and $\Sigma_3$ in the
following way: the conformal structure on $\Sigma \setminus \{z_1,
z_2, z_3\}$ can be described in terms of the {\it minimal area
metric} [Z] which we denote by $g_\Sigma$. This metric makes
$\Sigma$ as the union of three half cylinders $\Sigma_i$'s with
flat metric with each meridian circle has length $2 \pi$. The
metric is singular only at two points $p, \overline p \in \Sigma$
which lie on the boundary cycles of $\Sigma_i$. Therefore the
conformal structure induced from the minimal area metric naturally
extends over the two points $p, \overline p$. The resulting
conformal structure is nothing but the standard unique conformal
structure on $\Sigma= S^2 \setminus \{z_1, z_2, z_3\}$. We will
use this metric for the analytic estimates of pseudo-holomorphic
curves implicit in the argument. One important property of this
{\it singular} metric is that it is flat everywhere except at the
two points $p, \overline p$ where the metric is singular but
Lipschitz. Because the metric is Lipschitz at the two points and
smooth otherwise, it follows that the corresponding Sobolev
constants of the metric are all finite even though the metric is
singular at the two points $p, \, \overline p$ (see [\S 10, FOh]
for a similar consideration of such metric for the {\it open
string} version).

Following  [\S 10, FOh], [\S 11, En], we identify $\Sigma$ as the
union of $\Sigma_i$'s
$$
\Sigma = \cup_{i=1}^3 \Sigma_i \tag 4.3
$$
in the following way: if we identify $\Sigma_i$ with $(-\infty, 0]
\times S^1$, then there are 3 paths $\theta_i$ of length ${1 \over
2}$ for $i = 1, \, 2,\, \, 3$ in $\Sigma$ connecting $p$ to
$\overline p$ such that
$$
\aligned
\part_1\Sigma & = \theta_1\circ \theta_3^{-1} \\
\part_2 \Sigma & = \theta_2 \circ \theta_1^{-1} \\
\part_3 \Sigma & = \theta_3 \circ \theta_2^{-1}
\endaligned
\tag 4.4
$$
We fix a holomorphic identification of each $\Sigma_i, \, i= 1,2$
with $(-\infty, 0] \times S^1$ and $\Sigma_3$ with $[0,\infty)
\times S^1$ consider the decomposition (4.4).

We denote the identification by
$$
\varphi^+_i: \Sigma_i \to (-\infty, 0] \times S^1
$$
for the outgoing ends and
$$
\varphi^-_3: \Sigma_3 \to [0, \infty) \times S^1
$$
for the incoming end. We denote by $(\tau,t)$ the standard
cylindrical coordinates on the cylinders with identification $S^1
= \R/\Z$. So in our case, the length of the meridian has length 1
instead of $2\pi$.

We consider a cut-off function $\rho^-: (-\infty,0] \to [0,1]$ of
the type
$$
\rho^- = \cases 1 \quad & -\infty < \tau \leq -R_2 \\
0 \quad & -R_1\leq \tau \leq 0
\endcases
\tag 4.5
$$
and $\rho^+: [0,\infty) \to [0,1]$ by $\rho^+(\tau) =
\rho^-(-\tau)$ where $0 < R_1 < R_2 < \infty$ are arbitrary
numbers. We will just denote by $\rho$ these cut-off functions for
both cases when there is no danger of confusion.

We now consider the (topologically) trivial bundle $P \to \Sigma$ with fiber
isomorphic to $(M,\omega)$ and fix a trivialization
$$
\Phi_i: P_i:= P|_{\Sigma_i} \to \Sigma_i \times M
$$
on each $\Sigma_i$.
On each $P_i$, we consider the closed two form of the type
$$
\omega_{P_i}:= \Phi_i^*(\omega + d(\rho H_t dt)) \tag 4.6
$$
for a time periodic Hamiltonian $H: [0,1] \times M \to \R$. Note
that this naturally extends to a closed two form $\omega_P$
defined on the whole $P \mapsto \Sigma$ since the cut-off
functions vanish in the center region of $\Sigma$. Note that this
form is nondegenerate in the fiber and restricts to
$\Phi_i^*\omega$.

Such $\omega_P$ induces a canonical symplectic connection $\nabla
= \nabla_{\omega_P}$ [GLS], [En]. We consider the symplectic form
on $P$
$$
\Omega_{P,\lambda} := \omega_P + \lambda \omega_\Sigma \tag 4.7
$$
where $\omega_\Sigma$ is an area form and $\lambda > 0$ is a
sufficiently large constant. We will always normalize
$\omega_\Sigma$ so that $\int_\Sigma \omega_\Sigma = 1$ as before.

Next let $\widetilde J$ be the almost complex structure defined as
in (3.5) on each $\Sigma_i$ which naturally extends to the whole
$P$ due to the cut-off functions. This almost complex structure on
$P$ is $(H,J_0)$-compatible in that it satisfies

\roster
\item $\widetilde J$ is $\omega_P$-compatible on each fiber in that it
preserves the vertical tangent space
\item the projection $\pi: P \to \Sigma$ is pseudo-holomorphic, i.e,
$d\pi \circ \widetilde J = j \circ d\pi$.
\endroster
When we are given three $t$-periodic Hamiltonian $H =
(H_1,H_2,H_3)$ and a periodic family $J=(J^1,J^2,J^3)$,
$\widetilde J$ additionally satisfies

(3) For each $i$, $(\Phi_i)_*\widetilde J = j\oplus J_{H_i}$ where
$$
J_{H_i}(\tau,t,x) = (\phi^t_{H_i})^*J^i_t \tag 4.8
$$
on $M$ over the cylinder $\Sigma_i' \subset \Sigma_i$
in terms of the cylindrical coordinates. Here
$$
\Sigma_i' = \varphi_i^{-1}((-\infty, -R_{i,2}] \times S^1),
\quad i = 1,\, 2
$$
and
$$
\Sigma_3=\varphi^{-1}_3([R_{3,2}, \infty) \times S^1)
$$
for some $R_{i,2} > 0$. $\widetilde J$ also satisfies
$$
\widetilde J = j \oplus J_0
$$
on
$$
\Sigma_{cn} := \cup_{i=1}^3 C_i
$$
where $C_i \subset \Sigma_i$ is the
region corresponding to $[-R_{i,1},0]\times S^1$ for $i = 1,2$ and the one
corresponding to $[0, R_{3,1}]\times S^1$.
\medskip

The $\widetilde J$-holomorphic sections
$v$ over $\Sigma_i'$ are precisely the solutions of the equation
$$
{\part u \over \part \tau} + J^{i,\rho(\tau)}_t\Big({\part u \over
\part t} - X_{H_i}^{\rho(\tau)}(u)\Big) = 0 \tag 4.9
$$
if we write $v(\tau,t) = (\tau,t, u(\tau,t))$ in the trivialization
with respect to the cylindrical coordinates $(\tau,t)$ on $\Sigma_i'$
induced by $\phi_i^\pm$ above. In the center region $\Sigma_{cn}$,
they just become $J_0$ holomorphic curves
$$
u: \Sigma_{cn} \to M
$$
with respect to the given conformal structure $j$ on
$\Sigma_{cn}$.

Now we define the moduli space which will be relevant to the
definition of the pants product that we need to use. To simplify
the notations, we denote
$$
\widehat z =[z,w]
$$
in general and $\widehat z = (\widehat z_1, \widehat z_2, \widehat
z_3)$ where $\widehat z_i =[z_i,w_i] \in \text{Crit}\AA_{H_i}$ for
$i =1,2, 3$. We denote by 
$$
\widetilde \pi_2(M) = \text{Im}(\pi_2(M) \to H_2(M,\Z))
/\text{Tor}(H_2(M,\Z))
$$ 
the set of spherical homology classes (mod) torsion elements.

\definition{Definition 4.2} Consider the Hamiltonians
$H =\{H_i\}_{1\leq i \leq 3}$ and the closed two form $\omega_P$
on the trivial bundle $P = \Sigma \times M$ and its associated
connection $\nabla$. Let $\widetilde J$ be a $H$-compatible almost
complex structure.

\roster \item We denote by $\MM(H, \widetilde J; \widehat z)$ the
space of all $\widetilde J$-holomorphic sections $v: \Sigma \to P$
such that the closed surface obtained by capping off $u:=
pr_M\circ v$ with the discs $w_i$ taken with the same orientation
for $i = 1,2$ and the opposite one for $i =3$ represents zero in
$\widetilde \pi_2(M)$:
$$
[u \# (\bigcup_{i =1}^3 w_i)] = 0 \quad \text{in } \, \widetilde \pi_2(M). 
\tag 4.10
$$
Note that $\MM(H, \widetilde J; \widehat z)$ depends only on the
equivalence class of $\widehat z$'s:
\item
we say that $\widehat z'
\sim \widehat z$ if they satisfy
$$
z_i' = z_i, \quad w_i' = w_i \# A_i
$$
for $A_i \in \Gamma$ and $\sum_{i=1}^3 A_i = 0$ in $\widetilde
\pi_2(M)$.
\endroster
\enddefinition

The (virtual) dimension of $\MM(H, \widetilde J; \widehat z)$ is
given by
$$
\aligned \dim  \MM(H, \widetilde J; \widehat z) & = 2n -
(\mu_{H_1}(\widehat z_1) + n) - (\mu_{H_2}(\widehat z_2) +n)
- (-\mu_{H_3}(\widehat z_3) + n)\\
& = n + (\mu_{H_3}(\widehat z_3) - \mu_{H_1}(\widehat z_1) -
\mu_{H_2}(\widehat z_2)).
\endaligned
$$
Note that when $\dim \MM(H,\widetilde J;\widehat z) = 0$, we have
$$
n= -\mu_{H_3}(\widehat z_3) + \mu_{H_1}(\widehat z_1)
+ \mu_{H_2}(\widehat z_2)
$$
which is equivalent to
$$
k_3 = k_1 + k_2 - n
$$
if we write
$$
k_i = \mu_{H_i}(\widehat z_i). \tag 4.11
$$
This is exactly the grading of the Floer complex we adopted in
[Oh5]. Now the pair-of-pants product $*$ for the chains is defined
by
$$
\widehat z_1 * \widehat z_2 =  \sum_{\widehat z_3} \#(\MM(H,
\widetilde J;\widehat z)) \widehat z_3 \tag 4.12
$$
for the generators $\widehat z_i$ and then by linearly extending over
the chains in $CF_*(H_1) \otimes CF_*(H_2)$. Our grading convention makes
this product is of degree zero.

The following uniform energy bound will be used in our adiabatic
convergence argument in \S 6.

\proclaim{Proposition 4.3} Let $\widehat z_i=[z_i,w_i] \in
\text{Crit}\AA_{H_i}$ and $v: \Sigma \to M$ be a $\widetilde
J$-holomorphic section of $P \to \Sigma$ with the asymptotic
condition $\widehat z = (\widehat z_1,\widehat z_2,\widehat z_3)$.
Then we have the following inequality
$$
\align {1\over 2} \int_\Sigma |(Dv)^v|^2_{\widetilde J} 
& = \int v^*\omega_P -
\sum_{i=1}^2 \int_0^1 \rho_-'(\tau) H_i(t, u(\tau,t))\, dt -
\int_0^1 \rho'_+(\tau) H_3(t,u(\tau,t)) \tag 4.13 \\
& = - \AA_{H_3}([z_3,w_3]) + \AA_{H_1}([z_1,w_1]) + \AA_{H_2}
([z_2,w_2]) \\
& \quad - \sum_{i=1}^2 \int_0^1 \rho_-'(\tau) H_i(t, u(\tau,t))\,
dt - \int_0^1 \rho'_+(\tau) H_3(t,u(\tau,t)) \tag 4.14 \\
& \leq - \AA_{H_3}([z_3,w_3]) + \AA_{H_1}([z_1,w_1]) + \AA_{H_2}
([z_2,w_2]) \\
& \quad + \sum_{i=1}^2\int_0^1 \max H_{i,t}\, dt + \int_0^1 -\min
H_{3,t}\, dt \tag 4.15
\endalign
$$
In particular, we have the uniform upper bound for the vertical
energy of $v \in \MM(H,\widetilde J;\widehat z)$ which is
independent of the cut-off functions $\rho^\pm$ or of the choice
of $0 < R_1 < R_2 < \infty$.
\endproclaim
\demo{Proof} This is an immediate consequence of summing over $i =
1,2,3$ of (3.22) in Proposition 3.4. The only matter to be
clarified is the question on the orientation about the embedding
$$
\Sigma_i \hookrightarrow \Sigma.
$$
Recall that the embedding of the outgoing end
$$
\Sigma_i \cong (-\infty,0] \times S^1 \hookrightarrow \Sigma
\quad\text{for } \, i = 1,2
$$
is orientation preserving, but the embedding
$$
\Sigma_3 \cong [0,\infty) \times S^1 \hookrightarrow \Sigma
$$
is orientation reversing. The latter is responsible for the
negative sign in front of $\int \rho_+'(\tau) H_3(t,u(\tau,t))\,
dt$. The calculation leading to (3.22) was carried out for the
outgoing end $(-\infty,0] \times S^1$. Taking the orientation
change on the incoming end into consideration and summing (3.22)
over $i$, we get (4.13).

For the proof of (4.14), we just apply the identity from [\S 5,En]
$$
\int_\Sigma v^*\omega_P = -\AA_{H_3}([z_3,w_3]) +
\AA_{H_1}([z_1,w_1]) + \AA_{H_2} ([z_2,w_2]).
$$
(See also [Sc] but with different conventions).
This finishes the proof. \qed\enddemo

This general inequality can be also used the other way around.
More precisely, the following lower bound
will be a crucial ingredient for the limiting arguments
in the next section for the pseudo-holomorphic sections
with some asymptotic conditions which are allowed to vary
inside given Floer cycles.

\proclaim{Corollary 4.4} Let $\widehat z_i = [z_i,w_i]$ for $i=1,\, 2,\, 3$.
Suppose $\AA_{H_3}([z_3,w_3]) \geq c$ for some constant $c$. Then
we have the lower bound
$$
\AA_{H_1}([z_1,w_1]) + \AA_{H_2}([z_2,w_2]) \geq c -(E^+(H_1) +
E^+(H_2)) + E^-(H_3). \tag 4.16
$$
\endproclaim

\head{\bf \S 5. Construction of the $W^{1,2}$-sections: analysis
of the thick part}
\endhead

In this section, we will start with the proof of the main
existence result. We will treat the case $F=H$ in \S 8 in a more
careful way and improve the existence statement to prove [Theorem
3.11, Oh5] which was postponed to the present paper. We refer to
[Oh4,5] for more explantation on the Novikov cycles
$h_{\HH_1}(\alpha)$ or $h_{\HH_2}(\beta)$ here.

\proclaim{Theorem 5.1} Let $H,\,  F$ be as above and $J_0$ be any
compatible almost complex structure on $M$ and let $\delta, \,
\delta_1 > 0$ be given. Then for any $J' \in j_{(\phi,J_0)}$ there
exist some $[z,w] \in h_{\HH_1}(\alpha)$ and $[z', w'] \in
h_{\HH_2}(\beta)$ with
$$
\align
\AA_H([z,w]) & \leq \rho(H;1) + \frac{\delta}{2} \tag 5.1\\
\AA_{\widetilde F}([z',w']) & \leq \rho(\widetilde F;1) +
\frac{\delta}{2} \tag 5.2
\endalign
$$
for which the following alternative holds:

\roster
\item There exists a cusp solution $u: \R \times S^1 \to M$
$$
u = u^-\# u^+ : = (u^-_1\# u^-_2 \#\cdots u^-_{N_1})\# (u^+_1 \#
\cdots \# u^+_{N_2})
$$
where each $u^-_i$ is a cusp-curve consisting of a finite number
of  sphere bubbles and at most one principal component that
satisfies
$$
\dudtau + J  \Big(\dudt - X_F(u)\Big) = 0, \quad J^-_t=
(\phi_H^t)_*J'_t
$$
and  similarly for $u^+_j$'s whose principal components satisfy
$$
\dudtau + J  \Big(\dudt - X_H(u)\Big) = 0, \quad J^+_t=
(\phi_F^t)_*J'_t
$$
and $u^-_{N_1} \# u^+_1$ satisfies
$$
\cases \dudtau + J  \Big(\dudt - X_{(F,H)}(u)\Big) = 0 \\
u(\infty)=[\widetilde z',\widetilde w'], \quad u(-\infty) = [z,w]
\endcases
\tag 5.3
$$
on $(\R \times S^1) \setminus \{\tau = 0\}$
which lies in $W^{1,2}$ and satisfies
$$
\int (w^-)^*\omega + \int (u_-)^*\omega + \int (u_+)^*\omega -
\int (w^+)^*\omega = 0 \tag 5.4
$$
where $u^\pm$ are the parts on $\tau <0$ and $\tau > 0$
respectively, and $u$ is smooth near $(0,0)$ and satisfies
$$
\aligned & u^+_{N_2}(\infty)  = [\widetilde z',\widetilde w'],
\, u^-_1(-\infty) = [z,w],\\
& u^-_{N_1}(0,0)  = u^+_1(0,0) = q \in B(u)
\endaligned
\tag 5.5
$$
\item or
there is some $J'_t$-holomorphic sphere $v:S^2 \to M$ for some
$t\in [0,1]$ that passes through the point $q$ and in particular
$$
\aligned \AA_H([z,w]) + \AA_{\widetilde F}([z',  w']) & + \int_0^1
\max H_t\, dt + \int_0^1 -\min F_t\, dt \\
&\quad \geq A_S(\phi, J_0;J') - \delta_1.
\endaligned
\tag 5.6
$$
\endroster
\endproclaim

In this section and the next, we give the proof of Theorem 5.1.
Before launching the proof of Theorem 5.1, we would like to give
some heuristic explanation of our construction.  Morally we would
like to apply the pants product to the case
$$
H_1 = H, \, H_2 = \widetilde F, \, H_3 = 0 \tag 5.7
$$
and study the moduli space
$$
\MM(H, \widetilde J;\widehat z)
$$
of pseudo-holomorphic sections of an appropriate Hamiltonian
bundle $P \to \Sigma$, for the $\widetilde J$ chosen as in section 4,
that has the boundary condition
$$
\align u|_{\part_1\Sigma} & =[z,w] \in h_{\HH_1}(\alpha), \,
u|_{\part_2 \Sigma} =[z',w']  \in h_{\HH_2}(\beta)\\
u|_{\part_3 \Sigma} & = [q,\widehat q] \in \gamma
\endalign
$$
Here we denote $\widehat z = ([z,w],[z',w'];[q,\widehat q])$.
Note that because $H_3 = 0$ in the monodromy condition (5.8) and
the outgoing end with monodromy $H_2 = \widetilde F$ is equivalent
to the incoming end with monodromy $H_2 = F$, we can fill-up the
hole $z_3 \in \Sigma$ and consider the cylinder with one outgoing
and one incoming end with the monodromy $H$ and $F$ respectively.
In other words our Hamiltonian bundle $P \to \Sigma$ becomes just
the deformed mapping cylinder $E_f$ for
$$
f: \tau \to \phi^1_{\rho(\tau)}
$$
as defined in \S 2 where $\phi_s^t$ is the two parameter family
generated by the Hamiltonians
$$
s \in [0,1] \mapsto (1-s)H + sF
$$
for the limiting cut-off function
$$
\rho = \cases 0 \quad \text{for } \, \tau < 0 \\
1 \quad \text{for } \, \tau > 1.
\endcases
$$
This limiting moduli space is precisely the space of solutions of
(5.3). This heuristic reasoning is the motivation to the approach
that we take in the proof below.

\definition{Remark 5.2} We note that when $F = H$ and $\rho \equiv 1$
the Hamiltonian bundle $P$ is nothing but the mapping cylinder
$P = E_\phi$ with $\phi = \phi^1_H$ and the
the almost complex structure $\widetilde J$ is nothing but 
the pushforward of the complex structure 
$$
\widetilde J(\tau,t,x) = j\oplus J'(t,x)
$$
on $\R \times \R \times M$ under the covering projection
$$
\pi: \R \times \R \times M \to E_\phi 
=(\R \times \R \times M)/(\tau,t+1,\phi(x))
\sim (\tau,t,x).
$$
Furthermore, the natural connection induced by the Hamiltonian 
$H: [0,1] \times \R \to M$ is flat. This remark will be important
in section 7 for the proof of Theorem II.
\enddefinition

We now suppose that we are given any $J' \in j_{(\phi,J_0)}$, and one
fixed $J_3' \in j_{(\phi_{\e k}^1,J_0)}$. Let $(\phi,J_0)$ be as
before and $k$ be a generic Morse function. For the sake of
simplicity, we will also assume that $k$ has the unique minimum
point.

Instead of the triple (5.7) which is not allowed in defining the
pants product in the Floer theory because $H_3 = 0$ is degenerate,
we consider the admissible triple
$$
H_1 = H, \, H_2 = \widetilde{\e k \# F}, \, H_3 = \e k \tag 5.8
$$
for the nondegenerate $H$ and $F$. We note that if $F$ is
nondegenerate and $\e$ is sufficiently small, $(\e k)\# F$ is also
nondegenerate and there is a natural one-one correspondence
between $\text{Per}(H)$ and $\text{Per}((\e k)\# F)$ and the
critical points $\text{Crit}\AA_F$ and
$\text{Crit}\AA_{\widetilde{(\e k)\# F}}$. We will again denote by
$[z',w']$ the critical point of $\AA_{H_2}$ corresponding to
$[z',w'] \in \text{Crit}\AA_{\widetilde F}$.

Then we consider the triple of the families of almost complex
structures $J = (J_1, J_2, J_3)$ 
on the three ends $\Sigma_i\cong (-\infty,R] \times S^1$ 
under the trivialization
$$
\Phi_i: P|_{\Sigma_i} \to \Sigma_i \times M
$$
which are defined by
$$
\align J^1_t & = (\phi_H^t)_*J'_t\\
J^2_t & = (\phi^t_{\widetilde{\e k \# F}})_*(\phi^{1-t}_{\e
k})^*J'_{(1-t)}\\
J^3_t & = (\phi_{\e k}^t)_*J'_{3,t}
\endalign
$$
respectively.  Note that the family $(\phi^{1-t}_{\e k})^*J'_{(1-t)} \in
j_{((\phi^1_{\e k \# F})^{-1},J_0)}$ and $J_2$ is still a
$t$-periodic family. And as $\e \to 0$, $J_2$ converges to the
time reversal family of $t\mapsto (\phi^t_F)_*J'_t$
$$
t \mapsto (\phi^{(1-t)}_F)_*J'_{(1-t)}.
$$
We also choose $J'_{3,t}$ so that as $\e \to 0$, the path
converges to the constant path $J_0$ which is possible as
$\phi^1_{\e k}$ converges to the identity map.

With these choices, we study the limit as $\e \to 0$ of the moduli
space
$$
\MM(H^\e, \widetilde J^\e;\widehat z^\e)
$$
with the asymptotic boundary condition (5.4) for each $\e$. The
main difficulty in studying this limit problem in general is that
the limit is a {\it singular} limit and it is possible that the
image of pseudo-holomorphic curves can collapse to an object of
the Hausdorff dimension one in the limit (see [FOh1] for the study
of this kind of the singular limit problem). Our first non-trivial
task is to prove that there is a suitable limiting procedure for
any sequence $v^\e \in \MM(H^\e, \widetilde J^\e;\widehat z^\e)$
as $\e \to 0$. We will show that there exists a sub-sequence of
$v^\e \in \MM(H^\e, \widetilde J^\e;\widehat z^\e)$ that converges
to the union
$$
u \cup \chi
$$
in the Hausdorff topology, where $\chi$ is a negative (cusp)
gradient trajectory of $-k$ landing at the critical point $q \in
\text{Crit }f$ and $u$ is a solution of (6.20) but satisfying
$$
u(0,0) \in \text{Im } \chi
$$
instead of $u(0,0) = q$. In particular, we will prove that the
piece of Hausdorff dimension one is always a (negative) gradient
trajectory of $-k$.

This proven, it is easy to see that the above mentioned heuristic
discussion cannot produce a solution required in Theorem 5.1 {\it
unless the negative gradient trajectory of $-k$ landing at $q$ is
trivial} for some reason. At this stage, the condition
$$
\mu_{-\e f}^{Morse}(q) = 2n
$$
will play an essential role and enable us to conclude that the
only gradient trajectories of $-f$ {\it landing at} $q$ are the
constant map $q$ and so we are in a good position to start with.

Let $P \to \Sigma$, $\omega_P$ and $\Omega_{P,\lambda}$ be as in
\S 4. We equip $P$ with an $H$-compatible almost complex structure
$\widetilde J$ such that
$$
\widetilde J = j \oplus J'_i
$$
on each $\Sigma_i$ where
$$
J'_1 \in j_{(\phi,J_0)}, \quad J'_2=\widetilde J'_1 \in
j_{(\phi^{-1}, J_0)}, \quad J'_3 \in j_{(\phi_{\e f}^1,J_0)}. \tag
5.9
$$
More explicitly we define $\widetilde J$ by
$$
\widetilde J(\tau,t,x)(\alpha,\beta,\xi) = (-\beta, \alpha,
(\phi_{H_i^{\rho(\tau)}}^t)^*J^{\rho(\tau)}_{i,t}(\xi - \beta
X_{\rho(\tau)H_i}) + \alpha X_{\rho(\tau)H}) \tag 5.10
$$
on each $\Sigma_i$, where $s \in [0,1] \mapsto J_i^s$ is a fixed
path from $J_i^1 = \{(\phi_{H_i}^t)_*J_t'\}_{t\in [0,1]}$ to $J_0$
for each $i = 1,\, 2,\, 3$. Note that $\widetilde J$ is
$H$-compatible and naturally extends to the whole $P$ because of
the cut-off function $\rho$ and $J^s_i$ was chosen to be $J_i^0 =
J_0$ and the connection is trivial on $\Sigma_{cn}$.

The $\widetilde J$-holomorphic sections $v$ over $\Sigma_i$ are
precisely the solutions of the equation
$$
{\part u \over \part \tau} + J^{i,\rho(\tau)}_t\Big({\part u \over
\part t} - X_{H_i}^{\rho(\tau)}(u)\Big) = 0 \tag 5.11
$$
if we write $v(\tau,t) = (\tau,t, u(\tau,t))$ in the
trivialization with respect to the cylindrical coordinates
$(\tau,t)$ on $\Sigma_i'$ induced by $\phi_i^\pm$ above. In the
center region $\Sigma_{cn}$, they just become $J_0$-holomorphic
curves
$$
v: \Sigma_{cn} \to M
$$
with respect to the given conformal structure $j$ on
$\Sigma_{cn}$.

We denote $h = [\phi,H]$ and $f =[\phi, F]$ and let $\HH_1$ and
$\HH_2'$ be the paths from $\e k$ to $H$ and $\widetilde{(\e k)\#
F}$ respectively and $\HH_3$ the constant path $\e k$. Then by the
construction we have the identity
$$
1^\flat_{HF_*(\e k)} = h_{id}((1\cdot 1)^\flat) =
h_{\HH_1}(1^\flat)*h_{\HH'_2}(1^\flat)
$$
in the homology. In the level of cycles, for any cycles $\alpha
\in CF_*(H), \, \beta \in CF_* \widetilde{(\e k)\# F}$ and $\gamma
\in CF_*(\e k)$ with
$$
[\alpha] = 1^\flat, \, [\beta] = 1^\flat, \, [\gamma] = 1^\flat,
$$
we have
$$
h_{\HH_1}(\alpha)*h_{\HH_2}(\beta) = \gamma + \part_{\e k} (\eta)
\tag 5.12
$$
for some $\eta \in CF_*(\e k)$. We recall the following
Non-pushing down lemma

\proclaim{Lemma 5.3 [Lemma 6.8, Oh3]} Consider the cohomology
class $1 \in QH^*(M) = HQ^*(-\e k)$. Then there is the {\it
unique} Novikov cycle $\gamma$ of the form
$$
\gamma = \sum_j c_j [x_j,\widehat x_j] \tag 5.13
$$
with $x_j \in \text{Crit }(- \e k)$ that represents the class
$1^\flat = [M]$. Furthermore for any Novikov cycle $\beta$ with
$[\beta] = [M]$, i.e, $\beta$ with
$$
\beta = \gamma + \part_{\e k}( \delta )
$$
with $\delta$ a Novikov chain, we have
$$
\lambda_{\e k}(\beta) \geq \lambda_{\e k}(\gamma). \tag 5.14
$$
\endproclaim

It follows from Lemma 5.3 and (5.12) that
$$
\lambda_{\e k}(h_{\HH_1}(\alpha)*h_{\HH'_2}(\beta)) \geq -
{\delta\over 3} \tag 5.15
$$
with $\delta$ independent of $\e$
by choosing $\e$ sufficiently small. By the definition of

$\rho(\cdot; 1)$, we can find $\alpha$ and $\beta$ such that
$$
\aligned
 \rho(H;1) - {\delta \over 2} & \leq \lambda_{H}(h_{\HH_1}(\alpha))
\leq \rho(H;1) + {\delta \over 2} \\
\rho(\widetilde{\e k \# F};1) - {\delta \over 2} & \leq
\lambda_{\widetilde{\e k\# F}}(h_{\HH_2}(\beta))  \leq
\rho(\widetilde{\e k \# F};1) + {\delta \over 2}.
\endaligned
\tag 5.16
$$
Furthermore by the continuity of $\rho(\cdot;1)$, we also have
$$
\lim_{\e \to 0}\rho(\widetilde{\e k \# F};1)  = \rho(\widetilde
F;1). \tag 5.17
$$
Now we need to compare the levels of
$\lambda_{H}(h_{\HH_1}(\alpha))$ , $\lambda_{\widetilde{\e k\#
F}}(h_{\HH'_2}(\beta))$ and $\lambda_{\e
k}(h_{\HH_1}(\alpha)*h_{\HH'_2}(\beta))$ .

For each given $\e > 0$, let $[z_\e,w_\e]
\in h_{\HH_1}(\alpha)$  and $[z',w'] \in
h_{\HH'_2}(\beta)$ and $[q,\widehat q\# A] \in \gamma + \part_{\e
k}(\eta )$ such that the moduli space
$\MM(H^\e,\widetilde J^\e;\widehat z^\e)$ is non-empty. We may
choose $[q,\widehat q\# A] \in \gamma + \part_{\e k}(\eta)$ so
that
$$
\AA_{\e k}([q,\widehat q\# A]) \geq - {\delta \over 2} \tag 5.18
$$
using (5.15). Furthermore since
$[h_{\HH_1}(\alpha)*h_{\HH'_2}(\beta)] = 1_\flat$ and the {\it
unique} maximum point of $-\e k$ is homologically essential, we
may also assume that $q$ is the unique maximum point and $A=0$.
For each given $\e >0$,
let $v^\e$ be any $\widetilde J$-holomorphic section of $P \to
\Sigma$ that has the asymptotic boundary condition
$$
\aligned v|_{\part_1\Sigma} & =[z_\e,w_\e] \in h_{\HH_1}(\alpha), \,
v|_{\part_2 \Sigma} =[z'_\e,w'_\e]  \in h_{\HH'_2}(\beta)\\
v|_{\part_3 \Sigma} & = [q,\widehat q] \in \gamma.
\endaligned
\tag 5.19
$$
Then it follows from (4.16) of Corollary 4.4 that we have
$$
\AA_H([z_\e,w_\e]) + \AA_{\widetilde{\e k\# F}}([z'_\e,w'_\e])
\geq -\frac{\delta}{2} - ( E^+(H) + E^+(\widetilde{\e k\# F})) +
E^-(\e k).
$$
By choosing $\e > 0$ sufficiently small, we may assume that
$$
\AA_H([z_\e,w_\e]) + \AA_{\widetilde{\e k\# F}}([z'_\e,w'_\e])
\geq -(E^+(H) + E^+(\widetilde{\e \#F})) - \delta_2:=
C(H,\widetilde F) \tag 5.20
$$
where $\delta_2$ can be made arbitrarily small by letting $\e \to
0$. In particular, combining (5.16) and (5.20), we derive the
bounds
$$
\aligned
C(H,\widetilde F) - \rho(H;1) & \leq \AA_H([z_\e,w_\e]) \leq \rho(H;1) +
\frac{\delta}{2} \\
C(H,\widetilde F) - \rho(\widetilde{\e k\# F};1) &
\leq \AA_{\widetilde{\e\# F}}([z_\e,w_\e]) \leq \rho(\widetilde{\e k\# F};1) +
\frac{\delta}{2} \\
\endaligned
\tag 5.21
$$
for all $\e > 0$.  Therefore for each given $\e > 0$, the lower
bound of (5.21) implies that there are only finitely many possible
asymptotic periodic orbits among the generators $[z_\e,w_\e] \in
h_{\HH_1}(\alpha)$ and $[z'_\e,w'_\e] \in h_{\HH_2'}(\beta)$
respectively such that the corresponding moduli space
$$
\MM(H^\e, \widetilde J;\widehat z^\e)
$$
becomes non-empty for  the asymptotic condition (5.19).
On the other hand, since we assume that $F$ is generic nondegenerate
$\e k \# F$ are nondegenerate for all sufficiently small $\e >0$
and hence there are canonical one-one correspondences between
$\text{Crit}\AA_F$ and $\text{Crit}\AA_{\e k\# F}$. Therefore
as we let $\e \to 0$, we may assume that the asymptotic orbits
$[z'_\e, w'_\e] \in \text{Crit}\AA_{\e k\# F}$ converges to $[z',w']
\in \text{Crit}\AA_F$ and
$[z^\e,w^\e]\in \text{Crit}\AA_H$ to $[z,w] \in \text{Crit}\AA_H$
in $C^\infty$-topology. By taking the limit of (5.21) as $\e \to 0$,
we obtain (5.1) in Theorem 5.1.

Having this convergence statement made for the asymptotic orbits,
for the simplicity of notations, we will omit the subscript $\e$
from $[z_\e,w_\e]$ and $[z'_\e, w'_\e]$ and just denote them by
$[z,w]$ and $[z',w']$ from now on. We also note that because the
degree of $1$ is zero, we have
$$
\mu_H([z,w]) = \mu_{\widetilde{(\e k) \# F}}([z',w']) = \mu_{\e
k}([q,\widehat q]) = n
$$
and so
$$
\mu^{Morse}_{-\e k}(q) = n + \mu_{\e k}([q,\widehat q]) = 2n.
$$
In particular, any gradient trajectory $\chi: (K, \infty) \to M$
of $-\e k$ satisfying
$$
\dot\chi - \e \, \text{grad }k(\chi) = 0, \quad \lim_{\tau \to
\infty}\chi(\tau) = q
$$
must be the constant map $\chi \equiv q$, which is exactly what we
wanted to have in our heuristic discussion in the beginning of
this section.

We recall the uniform energy bound in Proposition 3.3 for the
vertical energy of the section $v$. Since this bound is uniform
for all sufficiently small $\e$ and $R_{i,2} > R_{i,1}$, we can
now carry out the adiabatic convergence argument for $\Sigma_3$.
To carry out this adiabatic convergence argument, we will
conformally change the metric on the base $\Sigma$ of the
fibration. We note that the vertical energy of the section $v$ is
invariant under the conformal change of the base metric. We will
realize this conformal change by a conformal diffeomorphism
$$
\psi_\e : C \setminus \{(0,0)\} \to (\Sigma, g_\e)
$$
where $C$ is the standard flat cylinder $\R \times S^1$ and the
metric $g_\e$ is constructed in a way similar to the minimal area
metric, but we will change the lengths of $\theta_i$'s in (4.4) in
the following way:
$$
\aligned
\text{length }\theta_1 & = 1 - {\e \over 2} \\
\text{length }\theta_2 & = \text{length }\theta_3 = {\e \over 2}.
\endaligned
\tag 5.22
$$
It is easy to see that we can choose the conformal diffeomorphism
$\psi_\e$ so that it restricts to a quasi-isometry
$$
\psi_\e: C \setminus D(\delta(\e)) \to (\Sigma_1 \cup \Sigma_2,
g_\e)
$$
and to a conformal diffeomorphism
$$
\psi_\e : D(\delta(\e)) \to (\Sigma_3, g_\e)
$$
where $D(\delta)$ is the disc around $(0,0)$ and $\delta(\e) \to
0$ as $\e \to 0$. We choose any sequence $\e_j \to 0$ and $v_j
\in \MM(H^{\e_j}, \widetilde J;\widehat z^{\e_j})$ that satisfy
(5.19).  Consider the compositions
$$
\widetilde v_j = v_j\circ \psi_{\e_j}; \, C \setminus \{(0,0)\}
\to P. \tag 5.23
$$
Since $\psi_\e$ is conformal, we have
$$
\int |(D\widetilde v_j)^v|^2_{\widetilde J} = \int |(D
v_j)^v|^2_{\widetilde J} \tag 5.24
$$
and hence the uniform energy bound
$$
\align {1 \over 2} \int |(D\widetilde v_j)^v|^2_{\widetilde J} &
\leq -\AA_{\e k}([q,\widehat q]) + \AA_{\widetilde{\e k \#
F}}([z',w']) +\AA_H([z,w]) \\
& \quad + \Big(\int_0^1 \max H \, dt + \int_0^1 \max \widetilde{\e
k\# F} \, dt + \int_0^1 - \min \e k\, dt\Big)\\
& \leq (\rho(H;1) + {\delta \over 2}) + (\rho(\widetilde{\e k\#
F}) + {\delta \over 2}) + {\delta \over 2} \\
& \quad + \Big(\int_0^1 \max H \, dt + \int_0^1 \max \widetilde{\e
k\# F} \, dt + \int_0^1 - \min \e k\, dt\Big) \\
& \leq \rho(H;1)+ \rho(\widetilde{F};1) + \int_0^1 \max H \, dt +
\int_0^1 \max \widetilde{ F} \, dt  + 2 \delta \tag 5.25
\endalign
$$
provided $\e$ is sufficiently small. Here we used (4.15) for the
first inequality, (5.16)-(5.18) for the second and the third
inequalities. There are two alternatives to consider: \roster
\item
there exists $c > 0$ such that
$$
{1 \over 2} \int_{D^2(\delta(\e_j))} |(D\widetilde
v_j)^v|^2_{\widetilde J} \geq c > 0 \tag 5.26
$$
for all sufficiently large $j$ or
\item
we have
$$
\limsup_{j \to \infty} {1 \over 2} \int_{D^2(\delta(\e_j))}
|(D\widetilde v_j)^v|^2_{\widetilde J} = 0. \tag 5.27
$$
\endroster

For the first alternative, $\widetilde v_j$ bubbles-off as $\e \to
0$. Since the bubble must be contained in a fiber by the maximum
principle by the $H$-compatibility of $\widetilde J$, it is
$J'_t$-holomorphic for some $t \in [0,1]$, we indeed have
$$
\limsup_{j \to \infty}{1 \over 2} \int |(D v_j)^v|^2_{\widetilde
J} = \limsup_{j \to \infty}{1 \over 2} \int |(D\widetilde
v_j)^v|^2_{\widetilde J} \geq A_S(\phi,J_0;J'). \tag 5.28
$$
Therefore we are in the second alternative in Theorem 5.1. It
remains to consider the case (5.27) which we will do in the next
section.

\head{\bf \S 6. Construction of the $W^{1,2}$-sections: analysis
of the thin parts} \endhead

We go back to the original sequence $v_j$. We can regard the
equation of pseudo-holomorphic sections with the asymptotic
conditions (5.19) as the following equivalent system of perturbed
Cauchy-Riemann equations:
$$
\align & \overline \part_{\widetilde J} v_i|_{C_i} = 0 \, \text{
for } \,
i = 1,2, 3 \tag 6.1-i \\
& v_1|_{\theta_3} = v_3|_{\theta_3},\, v_2|_{\theta_2} =
\widetilde v_1|_{\theta_2}, \, v_3|_{\theta_1} = v_2|_{\theta_1}
\tag 6.2
\endalign
$$
One can easily check that (5.19) are equivalent to the equation of
pseudo-holomorphic sections (6.1-i) with (6.2) because it is the
first order PDE of the Cauchy-Riemann type. We will study (6.1-i)
for each $i = 1,2,3$ separately with the matching condition (6.2)
in our mind.

We start with (6.1-3). It follows from (5.27) and by the conformal
invariance of the vertical energy that we have the energy bound
$$
{1 \over 2} \int_{\Sigma_3} |(D v_j)^v|^2_{\widetilde J} \leq
\delta_j \to 0 \tag 6.3
$$
as $j \to \infty$. On $[1, \infty) \times S^1$, (6.1-3) is
equivalent to
$$
\dudtau + J_t^\e \Big(\dudt - X_{\e f}\Big) = 0 \tag 6.4
$$
for the section $v(\tau,t) = (\tau,t, u(\tau,t))$ where $J_t^\e =
(\phi_{\e f}^t)_*J'_3$ with $J'_3 \in j_{(\phi^1_{\e k},J_0)}$.

We recall that $J^\e_t \to J_0$ as $\e \to 0$. The energy bound
(6.3) and the $\e$-regularity theorem (see [Corollary 3.4, Oh1]
for the context of pseudo-holomorphic curves), we immediately
derive the following uniform $C^1$-estimate.

\proclaim{Lemma 6.1} Let $\delta_j>0$ as in (6.1). Then there
exists $j_0$ such that for any $j \geq j_0$ and $u$, we have
$$
|Du_j(\tau, t)| \leq C \delta_j
$$
on $\Sigma_3 = \R_+ \times S^1$ with $C > 0$ independent of $j$.
In particular, we have
$$
\operatorname{length}(t\mapsto u_j(\tau,t)) \leq C \delta_j. \tag
6.5
$$
for any $\tau \in [0, \infty)$.
\endproclaim

We reparameterize $u_j$ to define
$$
\overline u_j(\tau,t) = u_j\Big({\tau \over \e_j}, {t\over
\e_j}\Big)
$$
on $[\e_j, \infty) \times \R /\e_j\Z$.  Then $\overline u_j$
satisfies
$$
{\part \overline u_j \over \part \tau}(\tau,t) + J_t^\e {\part
\overline u_j \over \part t}(\tau,t) - \operatorname{grad
}_{g_{J_t^\e}}f(\overline u_j)(\tau,t) = 0. \tag 6.6
$$
We will just denote $\operatorname{grad}$ for
$\operatorname{grad}_{g_{J_0}}$ from now on. We choose the `center
of mass' for each circle $t \mapsto u_j(\tau,t)$ which we denote
by $\chi_j(\tau)$. More precisely, $\chi_j(\tau)$ is defined by
the following standard lemma (see [K] for its proof).

\proclaim{Lemma 6.2} Suppose a circle $z: S^1 \to M$ has diameter
less than $\delta$ with $\delta= \delta(M)$ sufficiently small and
depending only on $M$. Then there exists a unique point $x_z$,
which we call the center of mass of $z$, such that
$$
z(t) = \exp_{x_z}\xi(t),\quad \int_{S^1}\xi(t)\, dt = 0. \tag 6.7
$$
Furthermore $x_z$ depends smoothly on $z$ but does not depend on
its parameterization.
\endproclaim

If $j_0$ is sufficiently large, (6.4) implies that the diameter of
all the circles $t \mapsto u_j(\tau,t)$  is less than $\delta$
given in  Lemma 6.2 for any $j$ and $\tau \in (0,\infty)$.
Therefore we can write $\overline u_j(\tau,t)$ as
$$
\overline u_j(\tau, t) = \exp_{\chi_j(\tau)} (\overline
\xi_j(\tau,t))
$$
with $\overline\xi_j(\tau,t) \in T_{\chi_j(\tau)}M$ for any $\tau
\in [1, \infty)$ for all $j \geq j_0$, $j_0$ sufficiently large.
Furthermore the $C^1$ estimate in Lemma 6.1 proves
$$
|\overline\xi(\tau,t)| \leq C \delta_j.
$$

We consider the exponential map
$$
\exp: \UU \subset TM \to M; \quad \exp(x, \xi):= \exp_x(\xi)
$$
and denote
$$
D_1\exp(x,\xi): T_xM \to T_xM
$$
the (covariant) partial derivative with respect to $x$ and
$$
d_2\exp(x, \xi): T_xM \to T_xM
$$
the usual derivative $d_2\exp(x,\xi):= T_{\xi}\exp_x: T_xM \to
T_xM$. We recall the property
$$
D_1\exp(x, 0) = d_2\exp(x, 0) = id
$$
which is easy to check. We now compute
$$
\aligned {\part \overline u\over \part \tau} & =
d_2\exp(\chi(\tau),\overline \xi(\tau,t))\Big({D \overline \xi
\over \part \tau}(\tau,t)\Big) + D_1\exp(\chi(\tau),
\overline \xi(\tau,t))(\dot\chi(\tau)) \\
{\part \overline u \over \part t} & = d_2\exp(\chi(\tau),
\overline \xi(\tau,t))\Big(
{\part \overline \xi \over \part t}(\tau,t)\Big) \\
\operatorname{grad} & f(\overline u)
=\operatorname{grad}f(\exp_{\chi(\tau)} \overline \xi(\tau,t)).
\endaligned
$$
Substituting these into (6.5) and multiplying by $
(d_2\exp(\chi_j(\tau), \overline \xi_j(\tau,t)))^{-1}$, we get
$$
\align {D \overline \xi_j \over \part \tau}(\tau,t) & +
(d_2\exp(\chi_j(\tau), \overline
\xi_j(\tau,t)))^{-1}(D_1\exp(\chi_j(\tau),\overline
\xi_j(\tau,t))(\dot\chi_j(\tau)))\\
& + (\exp_{\chi_j(\tau)})^*J_t^\e(\chi_j(\tau))\Big({\part
\overline \xi_j \over \part t}\Big)(\tau,t) -
(\exp_{\chi_j(\tau)})^*(\operatorname{grad}f)(\overline\xi_j(\tau,t))
= 0
\endalign
$$
on $T_{\chi_j(\tau)}M$. Using the center of mass condition (6.7)
and integrating over $t\in \R/\e_j\Z$,  we obtain
$$
\aligned \int_0^{\e_j} (d_2\exp(\chi_j(\tau), & \overline
\xi_j(\tau,t)))^{-1}
(D_1\exp(\chi_j(\tau),\overline \xi_j(\tau,t))(\dot\chi_j(\tau)))\, dt \\
& - \int_0^1 (\exp_{\chi_j(\tau)})^*(\operatorname{grad}f)
(\overline \xi_j(\tau,t)) \, dt= 0. \endaligned \tag 6.8
$$
Here we used the identities
$$
\align  & \int_{S^1} {\part \overline \xi_j \over \part t}(\tau,
t)\, dt = 0 \quad  \tag 6.9\\
& \int_{S^1} \overline \xi_j(\tau,t)\, dt \equiv 0 \equiv
\int_{S^1} {D \overline \xi_j \over \part \tau}(\tau, t)\, dt \tag
6.10
\endalign
$$
where the second identity of (6.10) is a consequence of the first.
On the other hand, it follows that there exists sufficiently small
$\delta_j > 0$ such that
$$
\align (1 - C\delta_j)|\dot\chi(\tau)| & \leq
|d_2\exp(\chi_j(\tau), \overline
\xi(\tau,t)))^{-1}(D_1\exp(\chi(\tau),\overline
\xi(\tau,t))(\dot\chi(\tau)))|\\
& \leq (1 + C\delta_j)|\dot\chi(\tau)| \tag 6.11
\endalign
$$
and
$$
|(\exp_{\chi(\tau)})^*(\operatorname{grad}f)(\xi(\tau,t))| \leq
(1+C\delta_j)|\operatorname{grad}f(\xi(\tau,t))| \leq (1+
C\delta_j)\|\operatorname{grad }f\|_{C^0}. \tag 6.12
$$
These follow from (6.7) and the following standard estimates on
the exponential map
$$
\aligned
|d_2\exp(x,\xi)(u)| & \leq C |\xi||u| \\
|D_1\exp(x,\xi)(u)| & \leq C |\xi||u|
\endaligned
$$
for $\xi,\, u \in T_xM$ where $C$ is independent of $\xi$, as long
as $|\xi|$ is sufficiently small, say smaller than injectivity
radius of the metric on $M$ (see [K]).

Combining (6.8), (6.11) and (6.12), we have obtained
$$
|\dot\chi_j(\tau)| \leq C\|\operatorname{grad }f\|_{C^0}.
$$
Therefore $\chi_j$ is equi-continuous and so on any given fixed
interval $[1, R] \subset [1, \infty)$, there exists a subsequence
of $\chi_j:[1,R] \to M$ uniformly converges  to some
$\chi_\infty:[1,R] \to M$. Furthermore it easily follows from
(6.12) and smoothness of the exponential map, we also have the
estimates
$$
\aligned |d_2\exp(\chi_j(\tau), \overline
\xi_j(\tau,t)))^{-1}(D_1\exp(\chi_j(\tau), & \overline
\xi_j(\tau,t))(\dot\chi_j(\tau)) - \dot\chi_j(\tau)| \\
& \leq C |\overline \xi_j(\tau,t)||\dot\chi_j(\tau)|
\endaligned
\tag 6.13
$$
and
$$
|\exp_{\chi_j(\tau)}^*(\operatorname{grad}f)(\overline\xi_j(\tau,t))
- \operatorname{grad}f(\chi_j(\tau))| \leq C
|\overline\xi_j(\tau,t)|. \tag 6.14
$$
where the constant $C$ depends only on $M$.  Since $\max|\overline
\xi_j|_{C^0} \leq C \delta_j \to 0$ uniformly, the equation (6.8)
converges uniformly to
$$
\dot\chi_\infty - \operatorname{grad }f(\chi_\infty) = 0.
$$
Therefore $\chi_\infty$ is a gradient trajectory of $f$. Recalling
the $C^1$-estimate in Lemma 6.1 and
$$
\overline u_j(\tau,t) = \exp_{\chi_j(\tau)} (\overline
\xi(\tau,t))
$$
we have proven that $\overline u_j|_{[0,R]\times \R/{\e_j\Z}}$
uniformly converges to $\chi_\infty$. By a boot-strap argument by
differentiating (6.8), this convergence can be turned into a
$C_\infty$-convergence.

Therefore by a standard argument of local convergence as in [Fl1],
the sequence $u_j$ converges to a {\it connected} finite union of
gradient trajectories of $-f$
$$
\chi_0 \# \chi_1 \#\cdots \# \chi_N
$$
for some $N \in \Z_+$ such that
$$
\chi_N(\infty) = q. \tag 6.15
$$
However since $\mu^{Morse}_{-\e f}(q) = 2n$, i.e, $q$ is a (local)
maximum of $-f$, $\chi_N$ must be constant which in turn implies
that all $\chi_j$ must be constant map $q$ for all $j$. Recalling
$$
u_j(\tau,t) = \overline u_j(\e_j\tau, \e_j t),
$$
we have proven that $u_j|_{\Sigma_3}$ converges to the constant
map $q$ in the $C^\infty$-topology.

Now we analyze the sequence $v_j$ over $\Sigma_1 \cup \Sigma_2$.
Consider a sequence $\e_j \to 0$ and $R_{3,2} = R_j \to 0$.
Obviously $R_{3,1}$ also goes to zero since $R_{3,2} > R_{3,1}$.
As we mentioned before the lower bound (5.21) implies 
by the definition of the Novikov chains that
there are only finitely many elements $[z,w] \in
h_{\HH_1}(\alpha)$ and $[z',w'] \in h_{\HH_2}(\beta)$ 
above the lower bound.
Therefore we can now apply the
Gromov-Floer type compactness arguments as $\e \to 0$.

When $v_j|_{\Sigma_1 \cup \Sigma_2}$ bubble-off, we are again in
the second alternative in Theorem 5.1. Therefore we assume that
there exists a constant $C > 0$ such that
$$
|Dv_j|_{\Sigma_1 \cup \Sigma_2} \leq C \tag 6.16
$$
by choosing a subsequence if necessary. Using this $C^1$ bound
together with the matching condition (6.2), we can bootstrap to
extract a (cusp)-limit $v_\infty$ on $\Sigma = \Sigma_1\cup \Sigma_2$
such that
$$
v_\infty = (v^-_1 \# v^-_2 \# \cdots \# v^-_N)\# (v^+_1 \# \cdots \# v^+_L)
$$
where $v^-_i$ for $1 \leq i \leq N-1$ and $v^+_j$ for $2 \leq i \leq L$  
are  maps from $\R \times [0,1]$ to $M$ 
but $v^-_N$ one from $(-\infty, 0] \times [0,1]$ and $v^+_1$ from 
$[0,\infty) \times [0,1]$. All of them satisfy 
$$
\cases {\part v \over \part \tau} + J'_t {\part v \over
\part t} = 0 \\
\phi(v(\tau,1)) = v(\tau,0), \quad \int |{\part v \over \part \tau
}|_{J'_t}^2 < \infty
\endcases
\tag 6.17
$$
In particular, the join $v^-_N \# v^+_1$ also satisfies (6.17) 
on $\R \times [0,1]$ with possible singularities along 
$\{\tau =0\}$. Here we would like to
point out that by the definition of the data 
$(P^\e, \widetilde J^\e, H^\e)$
and by the choice of minimal area metrics $g_\e$ in (5.22),
the bundle data $(P^\e, \widetilde J^\e)$ on $\R \times S^1 \setminus
\{(0,0)\}$ smoothly converges to $(P, \widetilde J)$ on compact sets
of $\R \times S^1 \setminus \{(0,0)\}$. Note that the latter datum
is smooth (see Remark 5.2). 
Now the matching condition (6.2) and the way how $v_\infty$
arises, especially the uniform $C^1$-bound (6.16) implies that
$v^-_N \# v^+_1$ becomes indeed smooth even across $\tau = 0$
by the ellipic regularity. And then since $v^-_N\# v^+_1$ has finite energy
on $\R \times \{(0,1)\}$ (as a $\widetilde J$-holomorphic section of $P$),  
it can be smoothly extended across $(0,0)$ by the removal singularity
theorem. Hence we have produced a smooth cusp-solution of 
the equation (6.17).

Now we consider $u_\infty: \R \times S^1 \to M$ defined by
$$
u_\infty(\tau,t) = \cases (\phi_H^t)(v_\infty(\tau,t)) \quad\text{for }\,
\tau <
0 \\
(\phi_F^t)(v_\infty(\tau,t)) \quad\text{for }\, \tau >0 \endcases
\tag 6.18
$$
The basic topological hypothesis
$$
[u \cup (\bigcup_{i =1}^3 w_i)] = 0 \quad \text{in } \, \widetilde\pi_2(M)
$$
in the definition of the pants product (see Definition 4.2) gives
rise to (5.4). Furthermore since $X_{(H,F)}$ smoothly matches near
$(0,0)$ by the basic assumption (2.1) and $u_\infty$ lies in $W^{1,2}$,
$u^-_\infty$ and $u^+_\infty$ smoothly matches around $(0,0)$ and $u(0,0) = q$.
Furthermore
$$
v(0,0) = u(0,0) = q \in B(u)
$$
while $v(\pm\infty) \in \text{Fix }\phi$ and hence $v$ cannot be
constant. This takes care of the alternative (1) in Proposition
5.2. This finally finishes the proof of Proposition 5.2 and so the
proof of Theorem 5.1.

\head{\bf \S 7. The case with $H_3 = H_1 \# H_2$}
\endhead

In this section, we specialize to the case $F=H$ in \S 6 and \S 6.
In this case, we can improve the arguments therein to prove the
following theorem, which is Theorem 3.11 in [Oh5] whose proof
was postponed to the present paper.

\proclaim{Theorem 7.1} Let $H$ and $J_0$ be as
before. And let $q \in \text{Int }B(u)$ and $\delta
>0$ be given. Then for any $J' \in j_{(\phi,J_0)}$, there exist
some generators $[z,w] \in h_{\HH}(\alpha)$ and $[z', w'] \in
h_{\HH}(\beta)$ with
$$
\aligned \AA_H([z,w]) & \leq \rho(H;1) + \frac{\delta}{2} \\
\AA_{\widetilde H}([z',w']) & \leq  \rho(\widetilde H;1) +
\frac{\delta}{2}
\endaligned
$$
such that the following alternative holds: \roster
\item The equation
$$
\cases \dudtau + J  \Big(\dudt - X_H(u)\Big) = 0 \\
u(\infty)=[\widetilde z',\widetilde w'], \quad u(-\infty) = [z,w]
\endcases
$$
has a cusp-solution
$$
u_1 \# u_2 \# \cdots \cdots \# u_N
$$
which is a connected union of Floer trajectories, possibly with
a finte number of sphere bubbles, for $H$ that
satisfies the conditions
$$
u_N(\infty) = [\widetilde z',\widetilde w'], \, u_1(-\infty) =
[z,w], \quad u_j(0,0) = q \in B(u).
$$
for some $1 \leq j \leq N$,
\item or there is some $J'_t$-holomorphic sphere $v:S^2 \to M$
for some $t\in [0,1]$ that passes through the point $q$.
\endroster
This in particular implies
$$
0 < A(\phi,J_0)\leq A(\phi)\leq \gamma(\phi) < \infty
$$
for any $\phi$ and $J_0$.
\endproclaim

The proof of this theorem will be identical except that we exploit
the fact that in the case where $F=H$,
we can actually construct a `flat' connection on $P
\to \Sigma$ by modifying the two parameter family
$$
\phi: [0,1] \times [0,1] \to \HH am(M,\omega)
$$
used in the construction of \S 3. We now explain this modification
here.

We start with $i = 2$ where $H_2 = \widetilde{\e k \# H}$. We will
interpolate the Hamiltonian isotopy
$$
\phi_2(0,t) = \phi_{\widetilde{\e k\# H}}^t
$$
and
$$
\phi_2(1,t) = \cases \phi_{\widetilde H}^{2t} \quad\text{for } \,
0
\leq t \leq {1 \over 2} \\
\phi^1_{\widetilde H}\cdot \phi_{\widetilde \e k}^{2t -1}
\quad\text{for } \, {1 \over 2} \leq t \leq 1\endcases
$$
by a family $\phi_s$ satisfying $\phi_s^1 = \phi_H^{-1}\cdot
\phi_{\e k}^{-1}$ for all $s \in [0,1]$. The corresponding family
is nothing but
$$
\phi_2(s,t) = \cases \phi_{\widetilde H}^{2t \over s}
\quad\text{for } \, 0
\leq t \leq {s \over 2} \\
\phi^1_{\widetilde H}\cdot \phi_{\widetilde{\e k}}^{2t -s \over
2-s} \quad\text{for } \, {s \over 2} \leq t \leq 1\endcases \tag
7.4
$$
For the case $i = 1$, we consider the family
$$
\phi_1(s,t) = \cases \phi_{\widetilde H}^{2t \over s}
\quad\text{for } \, 0
\leq t \leq {s \over 2} \\
\phi^1_{\widetilde H} \quad\text{for } \, {s \over 2} \leq t \leq
1\endcases \tag 7.5
$$
For the case $i = 3$, we define the family
$$
\phi_3(s,t) = \cases \phi_{\e k}^{2t \over s} \quad\text{for } \,
0 \leq t \leq {s \over 2} \\
\phi^1_{\e k} \quad\text{for } \, {s \over 2} \leq t \leq
1\endcases \tag 7.6
$$
Note that this family smoothly matches under the identification of
(4.4) except at the two points $p, \overline p$ where it is
continuous. These family define a flat connection on the
Hamiltonian fibration $P \mapsto \Sigma$ which is smooth away from
$p, \, \overline p$. For this family, it is straightforward to
compute the identity
$$
{1 \over 2} \int_\Sigma |(Dv)^v|^2 = -\AA_{\e k}([q,\widehat q]) +
\AA_{H}([z,w]) + \AA_{\widetilde{\e k \# H}} ([z',w']). \tag 7.7
$$
Once we have this, a simple examination of the proof in \S 6, 7
gives rise to the proof of Theorem 7.1.

\definition{Remark 7.2} In fact, the above construction of the flat
connection works for any triple $H = (H_1,H_2, H_3)$ satisfying
$$
H_3 = H_1\# H_2. \tag 7.8
$$
This fact was used by Schwarz [Sc] in his proof of the triangle
inequality. Since the details are not given in [Sc], we provide
the detail of this construction here.  For each $i$, we consider
the family
$$
\phi_1(s,t) = \cases \phi_{H_1}^{2t \over s} \quad\text{for } \, 0
\leq t \leq {s \over 2} \\
\phi^1_{H_1} \quad\text{for } \, {s \over 2} \leq t \leq
1\endcases \tag 7.9
$$
$$
\phi_2(s,t) = \cases \phi^1_{H_1} \quad\text{for } \, 0
\leq t \leq {s \over 2} \\
(\phi_H^1)\cdot\phi_{H_2}^{2t-s \over 2-s} \quad\text{for } \, {s
\over 2} \leq t \leq 1\endcases \tag 7.10
$$
$$
\phi_3(s,t) = \cases \phi_{H_1}^{2t \over s} \quad\text{for } \, 0
\leq t \leq {s \over 2} \\
\phi^1_{H_1} \cdot \phi^{2t -s \over 2-s}_{H_2} \quad\text{for }
\, {s \over 2} \leq t \leq 1\endcases \tag 7.11
$$
on $\Sigma_i$ for $i = 1,\, 2,\, 3$ respectively. It follows that
this smoothly matches under the identification of (4.4) away from
$p, \, \overline p$.
\enddefinition

\head {\bf References}
\endhead

\enddocument